\documentclass{article}
\usepackage[a4paper, margin=3cm]{geometry}
\usepackage[utf8]{inputenc}
\usepackage[english]{babel}

\usepackage{soul}
\usepackage{amsthm}
\usepackage{amssymb}
\usepackage{indentfirst}
\usepackage{amsmath}
\usepackage{graphicx, color}
\usepackage{caption}
\usepackage{enumitem}
\usepackage{algorithm}
\usepackage{algpseudocode}
\usepackage{authblk}
\usepackage{ulem}

\usepackage{tikz}

\graphicspath{{Imagenes/}}

\newtheorem{theorem}{Theorem}[section]
\newtheorem{lemma}[theorem]{Lemma}
\newtheorem{corollary}[theorem]{Corollary}
\newtheorem{proposition}[theorem]{Proposition}
\newtheorem{observation}[theorem]{Observation}
\newtheorem{definition}[theorem]{Definition}
\newtheorem{claim}{Claim}

\newcommand{\qitem}[1]{\noindent\leavevmode\hangindent1.5\parindent%
	\noindent\hbox to1.5\parindent{#1\hss}\ignorespaces}

\floatname{algorithm}{Procedure}

\newcommand{\exact}[2]{\hat{#1}^{[#2] }}
\newcommand{\stexact}[2]{\hat{#1}^{\{#2\} }}

\newcommand{\col}{{\rm col}}
\newcommand{\wcol}{{\rm wcol}}
\newcommand{\dcol}{{\rm dcol}}
\newcommand{\Reach}{{\rm Reach}}
\newcommand{\WR}{{\rm WReach}}
\newcommand{\DR}{{\rm DReach}}
\newcommand{\dr}{{\rm DR}_4}

\title{Colouring negative exact-distance graphs of signed graphs}

\author[1,2]{Reza Naserasr}
\author[3]{Patrice Ossona de Mendez} 
\author[4]{Daniel A. Quiroz}
\author[5]{Robert \v{S}\'amal} 
\author[1,*]{Weiqiang Yu}

\affil[1]{\small Department of Mathematics, Zhejiang Normal University, China.}
\affil[2]{\small Université Paris Cité, CNRS, IRIF, F-75013, Paris, France.}	
\affil[3]{\small Centre d'Analyse et de Math\'ematiques Sociales (CNRS, UMR 8557), Paris, France.}
\affil[4]{\small Instituto de Ingenier\'ia Matem\'atica-CIMFAV, Universidad de Valpara\'iso, Valpara\'iso, Chile.}
\affil[5]{\small Computer Science Institute of Charles University, Faculty of Mathematics and Physics, Charles University, Prague, Czech Republic.}
\affil[  ]{\footnotesize reza@irif.fr, pom@ehess.fr, daniel.quiroz@uv.cl, samal@iuuk.mff.cuni.cz, yuwq@zjnu.edu.cn.}
\affil[  ]{\footnotesize *Corresponding author.}

\date{}

\begin{document}
	
	\maketitle
	
	\begin{abstract}
		The \emph{$k$-th exact-distance graph},  of a graph $G$ has $V(G)$ as its vertex set, and $xy$ as an edge if and only if the distance between $x$ and $y$ is (exactly) $k$ in $G$. We consider two possible extensions of this notion for signed graphs. Finding the chromatic number of a negative exact-distance square of a signed graph is a weakening of the problem of finding the smallest target graph to which the signed graph has a sign-preserving homomorphism. We study the chromatic number of negative exact-distance graphs of signed graphs that are planar, and also the relation of these chromatic numbers with the generalised colouring numbers of the underlying graphs. Our results are related to a theorem of Alon and Marshall about homomorphisms of signed graphs.
		
	\end{abstract}

	\section{Introduction}

	The \underline{$k$-th power}, $G^k$, of a graph $G$ is the graph that has $V(G)$ as its vertex set, and $xy$ as an edge if and only if $d_G(x,y)\le k$. Problems related to the colouring of graph powers have received considerable attention, in part because of their connection to the frequency assignment problem in telecommunications. Of particular interest has been a conjecture of Wegner \cite{W77}, on the colouring of squares of planar graphs. Confirming a particular case of this conjecture, Thomassen \cite{Thomassen} proved that the square of every planar cubic graph is 7-colourable. Another landmark result on this topic is due to Agnarsson and Halld\'orsson \cite{AH03} who showed that there exists a constant $c_{d,k}$ such that for every $d$-degenerate graph we have $\chi (G^k)\le c_{d,k}\cdot\Delta(G)^{\lfloor k/2 \rfloor}$. Note that in this result, the exponent on $\Delta(G)$ is best possible as regular trees of radius $\lfloor k/2 \rfloor$ attest.
	
	The \underline{$k$-th exact-distance graph}, $G^{[\sharp k]}$, of a graph $G$ has $V(G)$ as its vertex set, and $xy$ as an edge if and only if $d_G(x,y)= k$. In recent years exact-distance graphs have also received considerable attention, especially after 
    the above-mentioned result of Agnarsson and Halld\'orsson has been refined, as we are going to describe next. 
    To state these new results formally, we need to define parameters that constitute a generalisation of the notion of degeneracy.
	
	Let $L$ be a total ordering of $V$, and $k\in\mathbb{N }\, \cup\{\infty\}$. We say that a vertex $x\in V$ is \underline{weakly $k$-reachable} from $y\in V$ if there exists an $xy$-path $P$ of length at most $k$ such that $x\le_{L} z$ for all vertices $z$ of $P$. If we additionally have $y\le_L z$ for all vertices $z\in P\setminus\{x\}$,  we say that $x$ is \underline{strongly $k$-reachable} from $y$. Let $\WR_k[G,L,y]$ and  $\Reach_k[G,L,y]$ be the sets of vertices that are weakly $k$-reachable and strongly $k$-reachable from $y$, respectively.  We set $$\wcol_k(G,L)=\max_{v\in V} |\WR_k[G,L,v]|,\phantom{space} \col_k(G,L)=\max_{v\in V} |\Reach_k[G,L,v]|,$$ and define the \underline{weak $k$-colouring number}, denoted $\wcol_{k}(G)$, and the \underline{strong $k$-colouring number}, denoted $\col_{k}(G)$, of a graph $G$, respectively, as follows:
	$$\wcol_{k}(G)=\min\limits_{L} \wcol_k(G,L), \phantom{space} \col_{k}(G)=\min\limits_{L} \col_k(G,L).$$
	
	These parameters, known as the generalised colouring numbers, are known to relate to other important graph parameters. For instance,  Grohe, Kreutzer, Rabinovich, Siebertz, and Stavropoulos  showed in \cite{Groheetal} that 
	\begin{equation}\label{eq:coltw}
		\col_1(G)\le \col_2(G)\le \dots \le \col_\infty (G)= \mathrm{tw}(G)+1,
	\end{equation}
	where $\mathrm{tw}(G)$ denotes the treewidth\footnote{A \underline{$k$-tree} is a graph which is either a clique of size $k + 1$ or is obtained from a smaller
		$k$-tree by adding a vertex adjacent to $k$ vertices which are pairwise adjacent. The \underline{treewidth}
		of a graph $G$ is the smallest $k$ such that $G$ is a subgraph of a $k$-tree.} of a graph. Also Zhu~\cite{Z09} showed that
	a graph class $\mathcal{C}$ has \underline{bounded expansion} if and only if for  every $k\in \mathbb{N}$, $\sup_{G\in \mathcal{C}}\wcol_{k}(G)<\infty$; classes with bounded expansion include planar graphs, classes excluding a fixed minor, and even those that exclude some fixed graph as a topological minor. 
	
	Ne\v{s}et\v{r}il and Ossona de Mendez \cite{NOdM12} proved that for every odd integer $k\ge 1$ and every class $\mathcal{C}$ with bounded expansion there is a constant $N_{k,\mathcal{C}}$ such that $\chi(G^{[\sharp k]})\le N_{k,\mathcal{C}}$ for every $G\in \mathcal{C}$. Van den Heuvel, Kierstead and Quiroz \cite{exactdistance} extended this result as follows.
	
	\begin{theorem}[Van den Heuvel, Kierstead and Quiroz~\cite{exactdistance}]\label{thm:unsignedexact}\mbox{}
		
		\qitem{a)}  For every odd positive integer $k$ and every graph $G$ we have $\chi(G^{[\sharp k]})\le \wcol_{2k-1}(G).$
		
		\qitem{b)}  For every even positive integer $k$ and every graph $G$ we have $\chi(G^{[\sharp k]})\le \wcol_{2k}(G)\cdot \Delta(G).$
	\end{theorem}

	
	Since, for instance, planar graphs have bounded expansion, this result implies that there is a constant $c$ such that every planar graph $G$ satisfies $\chi(G^k)\le c\cdot\Delta(G)^{\lfloor k/2 \rfloor}$, just as is implied by the result of Agnarsson and Halld\'orsson. But, indeed, we get a finer picture. As mentioned earlier, this result has contributed to exact-distance graphs receiving considerable attention in recent years \cite{bai, henry, Bousquet, Foucaud, la, priyam, quiroz}.

	
	The first goal of this paper is to present an extension of the theory using which one can replace part $(b)$ of Theorem~\ref{thm:unsignedexact} with a statement similar to that of part $(a)$, i.e., an upper bound independent of the maximum degree. As in many similar cases, this will be done using the notion of signed graphs.
 A \underline{signed graph} $(G, \sigma)$ is a pair, where $G$ is a graph and  $\sigma(e)\in\{+,-\}$ is the sign  of the edge $e$ in $G$. When the \underline{signature} $\sigma$ is of no particular relevance, we can denote the signed graph by $\hat{G}$. The signed graph on $G$ where all edges are negative is denoted by $(G, -)$.
 The sign of a path or a cycle in a signed graph is the product of the signs of all its edges, and the length of a shortest negative cycle is referred to be the negative girth.
 As in a considerable part of the literature, we interpret graphs as the subclass of signed graphs having all edges being negative. With this in mind, there seem to be two natural ways to extend the notion of exact-distance graphs to signed graphs, as follows.  
	The \underline{exact-distance $-k$ graph}, $\exact{G}{-k}$, of $\hat{G}$ is defined as the (unsigned) graph with $V(G)$ as its vertex set,
	and $xy$ as an edge if and only if $d(x, y)=k$ and every $xy$-path of length $k$ is negative in $\hat{G}$.  The \underline{strong exact-distance $-k$ graph}, $\stexact{G}{-k}$, of $\hat{G}$ has $V(G)$ as its vertex set,
	and $xy$ as an edge if and only if $d(x, y)=k$ and some $xy$-path of length $k$ is negative in $\hat{G}$. Clearly $E(\exact{G}{-k})\subseteq E(\stexact{G}{-k})$ and so $\chi(\exact{G}{-k})\leq \chi(\stexact{G}{-k})$. Moreover, if $\hat{G}$ represents a graph, that is, it has all edges negative, and $k$ is odd, then both $\exact{G}{-k}$ and $\stexact{G}{-k}$ are equal to the $k$-th exact-distance graph of the corresponding (unsigned) graph.
	
	The second goal is to see how the study of negative exact-distance graphs relates to, and sheds light into, the study of homomorphisms of signed graphs, particularly within planar graphs. A \emph{homomorphism} of a graph $G$ to a graph $H$ is a mapping of the vertices of $G$ to the vertices of $H$ such that adjacencies are preserved. When there is a homomorphism of $G$ to $H$, then we may write $G \to H$ or say $G$ \emph{maps} to $H$. Given a homomorphism $\varphi$ of $G$ to $H$, if $H$ has odd-girth at least $2k+1$, then the mapping $\varphi$ can be viewed as a proper colouring of $G^{[\sharp 2i+1]}$ for $i=1,\dots, k-1$ where vertices of $H$ are the colours. To capture the same for even values, we employ the notion of sign-preserving homomorphisms of signed graphs. That is a mapping of a signed graph $\hat{G}$ to a signed graph $\hat{H}$ that is not only a homomorphism of the underlying graphs but that also preserves the signs of the edges. It follows similarly that if $\phi$ is a homomorphism of $\hat{G}$ to $\hat{H}$ and $\hat{H}$ has negative girth at least $k$, then for all values of $i$, $1\leq i \leq k-1$, $\phi$ provides a colouring of $\stexact{G}{-i}$. The case of unsigned graphs is also captured by taking $(G,-)$ instead.

  
		 To apply these last observations to planar graphs we use the following theorem of Alon and Marshall.
		
		\begin{theorem}[Alon and Marshall \cite{AlonMarshall}]\label{thm:AM} Given an integer $k$, there exists a signed (simple) graph on $k2^{k-1}$ vertices which admits a sign-preserving homomorphism from any signed graph whose underlying graph admits an acyclic $k$-colouring\footnote{An \emph{acyclic $k$-colouring} of a graph $G$ is a proper $k$-colouring of its vertices in which every cycle receives at least three colours.}. 		
		\end{theorem}
		
		Since Borodin showed that every planar graph admits an acyclic 5-colouring \cite{borodin}, it follows that there is a signed (simple) graph on 80 vertices to which every signed planar graphs admits a sign-preserving homomorphism. Since in a sign-preserving homomorphism any two vertices of distance two which are joined by a negative path cannot map to the same vertex. Taking each vertex  of the target graph as a colour, we then obtain the following.
		 
	\begin{corollary}\label{coro:planar}
			For every planar signed graph $\hat{G}$ we have $\chi(\stexact{G}{-2})\leq 80.$
	\end{corollary}

 Kierstead and Yang \cite{KY} gave a short proof that $\col_2(G)$ is an upper bound for the acyclic chromatic number of $G$. Thus we also have the following.
 \begin{corollary}
 \label{coro:AM}
	For every signed graph  $\hat{G}$ we have
		$\chi(\stexact{G}{-2})\le \col_{2}(G)\cdot2^{\col_2(G)-1}$.     
 \end{corollary}

	\subsection{Our contributions}
	
	In this paper we study the chromatic number of exact-distance graphs of signed graphs $\hat{G}$ when~$G$ is planar. (Note that the exact-distance graph might not be planar any more: signed stars can produce arbitrarily large complete bipartite graphs.) We also study how the chromatic numbers of the exact-distance graphs of a signed graph $\hat{G}$ relate to the generalised colouring numbers of $G.$
 
 We start with a result for signed graphs of treewidth at most 2, that is, a subclass of planar graphs. Note that this result is for \textit{strong} negative exact-distance squares.
 
	\begin{theorem}\label{thm:treewidth2bound7}
		Let $\hat{G}$ be a signed graph. If $G$ has treewidth at most 2  then $\chi(\stexact{G}{-2})\le 7$.
	\end{theorem}

	
	Montejano, Ochem, Pinlou, Raspaud, and Sopena \cite{Montejanoetal} show that every signed graph of treewidth at most 2 admits a homomorphism to a 9-vertex signed graph. Moreover, they show that this is tight. Together with these results, Theorem~\ref{thm:treewidth2bound7} separates the notion of (strong) negative exact-distance graph and the notion of sign-preserving homomorphism, even within planar graphs.
	
	When proving Theorem \ref{thm:treewidth2bound7}  we actually prove that $\chi(G\cup\stexact{G}{-2})\le 7$, where $G\cup\stexact{G}{-2}$ is the graph with vertex set $V(G)$ and edge set $E(G)\cup E(\stexact{G}{-2})$. This result is tight; to see this let $G$ be the graph obtained by two negative paths $P_1$ and $P_2$ of length 2 and a universal vertex connected to all the vertices of $P_1$ through positive edges, and all those of $P_2$ through negative edges. It is easy to derive that $G$ has treewidth at most 2 and that $G\cup\stexact{G}{-2}$ is $K_7$, therefore $\chi(G\cup\stexact{G}{-2})=7$. 




	

 We now consider negative exact-distance graphs of general planar signed graphs. Noting that the upper bound of 80 given by Theorem~\ref{thm:AM} for planar graphs remains untouched despite many efforts, and that the upper bound in Corollary~\ref{coro:planar} applies also to $\chi(\exact{G}{-2})$, we improve this last bound from 80 to 77. More precisely, we prove the following.
 
	\begin{theorem}\label{thm:planarexact}
		For every signed planar graph $\hat{G}$ we have $\chi(\exact{G}{-2})\le 77$.
	\end{theorem}
	
	We remark that the proof of Theorem \ref{thm:planarexact} is essentially self-contained and, in particular, does not rely on Borodin's result. Our techniques can be further used to obtain a bound of 76, but for simplicity of presentation we abstain from doing so here.
	


	So far we have only considered distance 2, and now move to larger distances. If $k\ge 4$ is even we cannot have constant upper bounds on $\chi(\stexact{G}{-k})$, even for the class of outerplanar graphs.
	To see this, let $\hat{G}$ be obtained from a star $K_{1,\ell}$ by replacing every edge with both a positive path of length $k/2$ and a negative path of length $k/2$. (This is not possible when $k=2$, as we assume graphs are simple.) Then clearly, $\stexact{G}{-k}$ contains a clique $K_\ell$, which implies that $\chi(\stexact{G}{-k})\ge \ell$. This tells us that we cannot obtain upper bounds like the one of Corollary~\ref{coro:AM} for $\chi(\stexact{G}{-k})$ when $k\ge 4$ is even. However we do give similar upper bounds for $\chi(\exact{G}{-k})$ in the two following results. Firstly, we generalise part $(a)$ of Theorem~\ref{thm:unsignedexact}. This result is our starting point towards proving Theorem \ref{thm:planarexact}.
	
	\begin{theorem}\label{relation}
		For every signed graph $\hat{G}$ and positive integer $k$ we have
		$\chi(\exact{G}{-k})\le \wcol_{2k}(G)$. Moreover if $k$ is odd, we have $\chi(\exact{G}{-k})\le \wcol_{2k-1}(G)$.
	\end{theorem}

Corollary \ref{coro:AM} tells us that there is a function $f$ such that $\chi(\exact{G}{-2})\le f(\wcol_2(G))$. Thus it is natural to hope that $\chi(\exact{G}{-k})$ could be bounded in terms of $\wcol_k(G)$ and not just in terms of $\wcol_{2k}(G)$, as we get from Theorem~\ref{relation}. We show that this is indeed the case. 

\begin{theorem}\label{functionwcolk}
		For every signed graph $\hat{G}$ and positive integer $k$ we have
		$$
          \chi(\exact{G}{-k})\le (2\lceil k/2\rceil+2)^{\wcol_{k}(G)}.
        $$ 

\end{theorem}


	This upper bound in terms of $\wcol_{k}(G)$ is tight in the sense that it is not possible to give such a bound depending only on $\wcol_{\ell}(G)$ such that $\ell<k$. To see this, for $n,k\geq 2$, let $\hat{S}_{n,k}$ be the signed graph obtained from $K_n$ by replacing each edge of $K_n$ with a negative path of length $k$ (and $S_{n,k}$ its underlying graph). Then obviously, $\chi(\hat{S}_{n,k}^{[-k]})=n$. However, we have $\wcol_{k-1}(S_{n,k})\leq k+1$. To verify this, we order the vertices in the following way: first order the original vertices, those of $K_n$, arbitrarily, and later the added degree-$2$ vertices of the negative paths. Clearly, for each original vertex~$v$, there is no vertex that is weakly $(k-1)$-reachable from~$v$. And each added degree-$2$ vertex can only weakly $(k-1)$-reach the vertices on the same negative path, which implies that $\wcol_{k-1}(S_{n,k})\leq k+1$. 

		
		

	The rest of the paper is organized as follows. In Section \ref{sec:treewidth} we prove Theorem \ref{thm:treewidth2bound7}. In Section~\ref{sec:largedistance} we prove Theorem~\ref{relation} and Theorem~\ref{functionwcolk}. We prove Theorem~\ref{thm:planarexact} in Section~\ref{sec:planar}. We conclude the paper in Section \ref{sec:remarks} with some remarks about signed graphs with large treewidth. We end the introduction with some notation.

	For a positive integer $k$,
	let $[k]=\{1,2,\cdots,k\}$.
	For a vertex $v\in V$,
	denote by $N^{k}(v)$ the \underline{$k$-th neighbourhood} of $v$,
	that is, the set of vertices different from $v$ with distance at most $k$ from~$v$. We also set $N^{k}[v]=N^{k}(v)\cup \{v\}$. For a path $P$ and vertices $x, y\in V(P)$, we denote by $xPy$ the subpath of $P$ with endpoints $x$ and $y$. 
	Let $P$ and $P'$ be two paths with $x,y\in V(P)$, $u,v\in V(P')$ and $yu\in E(G)$, then $xPy$-$uP'v$ is the walk obtained by concatenating the paths $xPy$ and $uP'v$. For a walk $W$ we use $||W||$ to denote its length. 
	

	\section{Strong exact-distance squares in graphs of treewidth at most 2}\label{sec:treewidth}

	
	In this section we prove  Theorem~\ref{thm:treewidth2bound7}.
	As mentioned in the introduction, we actually prove that $\chi(G\cup\stexact{G}{-2})\le 7$. Notice that the deletion of an edge in $\hat{G}$ cannot add any edge to $G\cup\stexact{G}{-2}$, so we may assume that $G$ is a 2-tree. 
	
 
	We will construct an assignment $\mathcal{C}$ on $V(G)$ with $\mathcal{C}(x)=(c(x), A_x, B_x)$, where $c(x)\in [7]$ and $A_x, B_x\subset [7]$, such that
	\begin{itemize}
		\item[(i)]  $c(x)\notin A_x\cup B_x$, 
		\item[(ii)] $A_x\cap B_x= \varnothing$,
		\item[(iii)] $|A_x|=|B_x|=3$.
	\end{itemize}
	Moreover, if $xz \in E(G)$ the assignment $\mathcal{C}$ will satisfy  each of the following properties
	\begin{itemize}
		\item[(iv)] $A_x\cap A_z$, $A_x\cap B_z$, $B_x\cap A_z$, $B_x\cap B_z$ are all non empty, 
		\item[(v)] $c(x)\in A_z$, $c(z)\in A_x$ if $xz$ is positive,
		\item[(vi)] $c(x)\in B_z$, $c(z)\in B_x$ if $xz$ is negative.
	\end{itemize}
	The first entry of $\mathcal{C}(x)$, namely $c(x)$ will be the colour of $x$. Intuitively then, by properties (v) and (vi), $A_x$ represents the colours available for the positive neighbours of $x$, while $B_x$ represents those available for the negative neighbours of $x$.

	Before showing that such an assignment is possible, we show that the function $c\colon V(G) \rightarrow \{1, \dots , 7\}$ obtained in this assignment is a proper colouring of $G\cup\stexact{G}{-2}$. For each edge $xy\in E(G\cup\stexact{G}{-2})$, either $xy\in E(G)$, or $x$ and $y$ have a common neighbour $z$ such that one of the edges $xz$ and $yz$ is positive, and the other is negative. Say we are in this latter case, with $xz$ being positive and $yz$ being negative, then by (v) $c(x)\in A_z$ and by (vi) $c(y)\in B_z$, which together with (ii) implies that  $c(x)\ne c(y)$. Otherwise, $xy\in E(G)$, by conditions (i), (v), and (vi), we obtain that $c(x)\ne c(y)$, showing that we would have a proper colouring of $G\cup\stexact{G}{-2}$.

	Now we show by induction on $|V(G)|$ that such an assignment can be obtained. For the base case $|V(G)|=2$, say $V(G)=\{x,y\}$. If $\sigma(xy)=+$, then $\mathcal{C}(x)=(1,\{2,3,4\},\{5,6,7\})$ and $\mathcal{C}(y)=(2,\{1,3,5\},\{4,6,7\})$ is an assignment satisfying the conditions. Otherwise, the assignment $\mathcal{C}(x)=(1,\{2,3,4\},\{5,6,7\})$ and $\mathcal{C}(y)=(5,\{2,4,6\},\{1,3,7\})$ satisfies the conditions.
	
	Now we take a signed graph $(G,\sigma)$ of $|V(G)|\geq 3$. Since $G$ is a 2-tree, we can choose a vertex $z$ such that $d(z)=2$ and let $N(z)=\{x,y\}$. By induction we know that $G-z$ admits an assignment satisfying conditions (i)--(vi), we now extend the assignment to the whole graph. Since $G$ is a $2$-tree, $x$ and $y$ are adjacent. Without loss of generality, let $\mathcal{C}(x)=(x_1,\{x_2, x_3, x_4\}, \{x_5, x_6, x_7\})$ and $\mathcal{C}(y)=(x_2,\{x_1,x_3,x_5\},\{x_4,x_6,x_7\})$ when $\sigma(xy)=+$; let $\mathcal{C}(x)=(x_1,\{x_2, x_3, x_4\}, \{x_5, x_6, x_7\})$ and $\mathcal{C}(y)=(x_5,\{x_2,x_6,x_4\},\{x_1,x_3,x_7\})$ when $\sigma(xy)=-$. We then consider subcases depending on the signs of the edges $xz$ and $yz$. It is not hard to check that Table \ref{table:signs} gives, for each subcase, an assignment $\mathcal{C}(z)$ that satisfies the conditions (i)--(vi) in relation to $\mathcal{C}(x)$ and $\mathcal{C}(y)$, extending $\mathcal{C}$ in the desired way.
	\begin{table}[H]
	\begin{center}
		\begin{tabular}{|c|c|}
			\hline
			$(\sigma(xy), \sigma(xz),\sigma(yz))$&	$\mathcal{C}(z)$ \\
			\hline
			$(+,+,+)$& $(x_3,\{x_1,x_2,x_6\}, \{x_4,x_5,x_7\})$\\
			
			$(+,+,-)$& $(x_4,\{x_1,x_3,x_6\}, \{x_2,x_5,x_7\})$\\
			
			$(+,-,+)$& $(x_5,\{x_2,x_3,x_7\}, \{x_1,x_4,x_6\})$\\
			
			$(+,-,-)$& $(x_7,\{x_3,x_4,x_5\}, \{x_1,x_2,x_6\})$\\
			\hline
			$(-,+,+)$& $(x_2,\{x_1,x_4,x_5\}, \{x_3,x_6,x_7\})$\\
			
			$(-,+,-)$& $(x_3,\{x_1,x_2,x_6\}, \{x_4,x_5,x_7\})$\\
			
			$(-,-,+)$& $(x_6,\{x_2,x_3,x_5\}, \{x_1,x_4,x_7\})$\\
			
			$(-,-,-)$& $(x_7,\{x_3,x_4,x_6\}, \{x_1,x_2,x_5\})$\\
			\hline
			
		\end{tabular}
  \caption{Assignment $\mathcal{C}(z)$ for each case}
  \label{table:signs}
	\end{center}
	\end{table}
	
	Thus we are able to extend the assignment $\mathcal{C}$ to $z$ in every case, and the result follows.
	
	In this proof, we have in fact done something stronger as explained next without giving detailed proof. As we will see, we have built a signed graph on $7\times \binom{6}{3}=140$ vertices with the following two properties: the first is that it admits a sign-preserving homomorphism from any signed 2-tree and the second is that the target graph itself can be coloured with seven colours in such way that adjacent vertices or vertices connected with a negative path of length 2 are assigned distinct colours.
	
	The target graph, denoted $\hat{P}_{1,3,3}$, has as its vertices all ordered partitions of  $\{1,2,3,4,5,6,7\}$ into three sets first of order 1, second and third both of order 3. (Thus the first entry has 7 options, each of them with $\binom{6}{3}$ possibilities, therefore, in total there are 140 vertices.) A partition $(x, A, B)$ is adjacent to $(y, A', B')$ with a positive edge if $x\in A'$ and $y\in A$. They are adjacent  with a negative edge if $x\in B'$ and $y\in B$. The proof given above can be read as claiming that any signed 2-tree admits a homomorphism to $\hat{P}_{1,3,3}$. To complete the claim one can check that the assignment of $x$ to vertex $(x,A, B)$ is a 7-colouring admitting all the required properties.



	\section{Larger distances and the generalised colouring numbers}\label{sec:largedistance}
	
	We first prove Theorem \ref{relation}, in fact a stronger version of it, and for this we introduce a refined version of the weak colouring numbers introduced by Van den Heuvel, Kierstead and Quiroz \cite{exactdistance}.

	Let $G=(V,E)$ be a graph,
	$L$ a total ordering of $V$,
	and $k$ a positive integer.
	For a vertex $y\in V$,
	let $\DR_k[G,L,y]$ be the set of vertices $x$ such that there is an $xy$-path $P_{x}=z_{0},\cdots,z_{s}$,
	with $x=z_{0}$, $y=z_{s}$,
	of length $s\le k$,
	such that $x$ is the minimum vertex in $P_{x}$ with respect to $L$,
	and such that $y\le _{L} z_{i}$ for $\lfloor \frac{1}{2}k \rfloor +1\le i \le s$.
	The distance-$k$-colouring number $\dcol_{k}(G)$ of a graph $G$ is defined as follows:
	
	$$\dcol_{k}(G,L)=\max_{v\in V} |\DR_k[G,L,v]|,$$
	
	$$\dcol_{k}(G)=\min\limits_{L} \dcol_k(G,L).$$
	
	
	Since $\Reach_k[G,L,y]\subseteq \DR_k[G,L,y]\subseteq \WR_k[G,L,y]$ for every ordering $L$,
	distance $k$ and vertex $y$,
	we have that $\col_{k}(G)\le \dcol_{k}(G)\le \wcol_{k}(G)$.
	In order to prove Theorem~\ref{relation},
	we actually prove the following stronger theorem.

	\begin{theorem}\label{stronger}
		For a signed graph $(G, \sigma)$ and a positive integer $k$ we have
		$\chi(\exact{G}{-k})\le \dcol_{2k}(G)$. Moreover if $k$ is odd, we have $\chi(\exact{G}{-k})\le \dcol_{2k-1}(G)$.
	\end{theorem}
	
	\proof Our proof builds on the proof of Van den Heuvel, Kierstead and Quiroz \cite[Theorem 2.1]{exactdistance}. We prove that $\chi(\exact{G}{-k})\le \dcol_{2k}(G)$, and leave the improvement for odd $k$ to the reader. For a positive integer $k$ and graph $G=(V,E)$,
	set $p=\dcol_{2k}(G)$ and let $L$ be an ordering of $V$ that witnesses $\dcol_{2k}(G,L)=\dcol_{2k}(G)$. 
	Moving along the ordering $L$ we assign to each vertex $y\in V$ a colour $a(y)\in [p]$ that is different from $a(x)$ for all $x\in \DR_{2k}[G,L,y]\setminus \{y\}$. Let $A^{\lfloor k/2\rfloor}[y]= N^{\lfloor k/2\rfloor}[y]$ if $k$ is odd, and
	$A^{\lfloor k/2\rfloor}[y]= N^{\lfloor (k-1)/2 \rfloor}[y]\cup \{v\in V\mid d(v,y)=k/2$ and every $vy$-path of length $k/2$ is positive$\}$, when $k$ is even.
	Define $\mu(y)$ as the minimum vertex with respect to $L$  in $A^{\lfloor k/2\rfloor}[y]$.
	Then define $c: V\rightarrow [p]$ by $c(y)=a(\mu(y))$.
	We claim that $c$ is a proper $p$-colouring of $\exact{G}{-k}$.
	
	Consider any edge $e=uv$ in $\exact{G}{-k}$.
	Then there is a negative path of length $k$,
	$P=x_{0},\cdots,x_{k}$ with $x_{0}=u$ and $x_{k}=v$.
	We first show that $\mu(u)\neq \mu(v)$. If $k$ is odd then $A^{\lfloor k/2\rfloor}[u]\cap A^{\lfloor k/2\rfloor}[v]=\varnothing$ so we indeed have $\mu(u)\neq \mu(v)$.
	If $k$ is even and we have $\mu(u)=\mu(v)$ ,
	then we also have $d(u,\mu(u))=d(v,\mu(u))=k/2$. Moreover, by definition of $\mu(u)$, there is a positive $u\mu(u)$-path and a positive $v\mu(v)$-path, both of length $k/2$. However together these two paths form a positive $uv$-path of length $k$, contradicting the fact that $uv$ is an edge of $\exact{G}{-k}$.
	
	
	%
	
	To complete the proof, i.e., to show that $c(u)\ne c(v)$, it remains to show that $a(\mu(u))\ne a(\mu(v))$. Since $\mu(u), x_{\lfloor k/2 \rfloor}\in N^{\lfloor k/2\rfloor}[u]$,
	there exist a path $P_1$ between $\mu(u)$ and $x_{\lfloor k/2\rfloor}$ of length at most $k$ such that $V(P_1)\subseteq N^{\lfloor k/2\rfloor}[u]$. 
	Similarly, there exist a path $P_2$ between $x_{\lceil k/2 \rceil}$ and $\mu(v)$ of length at most $k$ such that $V(P_2)\subseteq N^{\lfloor k/2 \rfloor}[v]$.

	Then the union of $P_1$ and $P_2$ (together with the edge $x_{\lfloor k/2 \rfloor}x_{\lceil k/2 \rceil}$ if $k$ is odd) forms a path $P^*$ between $\mu(u)$ and $\mu(v)$ of length at most $2k$. (For the improvement when $k$ is odd, note that in this case the length is at most $2k-1$.)
	
	We may assume without loss of generality $\mu(u)<_L \mu(v)$.
	We claim that we have $\mu(u)\in \DR_{2k}[G,L,\mu(v)]$, which implies $a(\mu(u))\ne a(\mu(v))$, as desired. The witness for this reachability is $P^*$.
    Observe that all the vertices of $P_1$ but possibly $x_{\lceil k/2 \rceil}$ are in $A^{\lfloor k/2\rfloor}[u]$ and that $\mu(u)$ is the minimum of this set  with respect to the ordering $L$. In this case that $x_{\lceil k/2 \rceil}$ is not in $A^{\lfloor k/2\rfloor}[u]$ there must be a negative path of length $ \lceil k/2 \rceil $ connecting $u$ to $x_{\lceil k/2 \rceil}$. Similarly, with respect to the ordering $L$, $\mu(v)$ is smaller than all vertices of $P_2$ but possibly $x_{\lceil k/2 \rceil}$ in the special case.
    Therefore, all the vertices of $P^*$ satisfy the conditions that $\mu(u)\in \DR_{2k}[G,L,\mu(v)]$ except for the possibility of $x_{\lceil k/2 \rceil}$ being smaller than $\mu(u)$ with respect to $L$ in the case when $k$ is even. But if that is the case, then $x_{\lceil k/2 \rceil}$ is also smaller than $\mu(v)$, and then there are negative paths of length $k/2$ between this vertex and each of $u$ and $v$, implying a positive path of length $k$ between $u$ and $v$, contradicting the fact that $uv\in \exact{G}{-k}$.
	\qed


	\subsection{A bound through shorter reachability}
	
	In this section we prove Theorem \ref{functionwcolk}.

    Let $G=(V,E)$, we take $L$ as an ordering of $V=\{x_1,\dots , x_n\}$ that witnesses $\wcol_{k}(G,L)=\wcol_{k}(G)$. Moving along the ordering $L$, we greedily colour $V$ with $f$: $V\rightarrow [\wcol_{k}(G)]$, that is, we give to each vertex $v$ a colour $f(v)$ that is different to $f(u)$ whenever $u\in \WR_{k}[G,L,v]\setminus\{v\}$. For a vertex $u\in \WR_{\lfloor k/2 \rfloor}[G,L,v]$, we let $sdist(u,v)$ be their signed distance, where $|sdist(u,v)|$ equals to the length of the shortest path that witnesses $u\in \WR_{\lfloor k/2 \rfloor}[G,L,v]$, and its sign being positive if some shortest paths are positive and negative if all shortest paths are negative. To each vertex $v$ we assign a colour $\alpha(v)$, which is a vector of length $\wcol_{k}(G)$, where the $i$-th entry is $sdist(u,v)$ if there exists a vertex $u\in \WR_{\lfloor k/2 \rfloor}[G,L,v]$ and $f(u)=i$, the entries not defined are set to $\ast$. Note that the number of colours used is at most $(2\lceil k/2\rceil+2)^{\wcol_{k}(G)}$. So it suffices to show that $\alpha$ is a proper colouring of $\exact{G}{-k}$.

	Consider any edge $uv$ in $\exact{G}{-k}$. By definition there is a negative $uv$-path $P$ of length $k$ in $G$. Let $w$ be the minimum vertex in $P$ with respect to $L$. We assume without loss of generality that $d_{P}(w, u)\leq d_{P}(w, v)$, then $P(w,u)$ is the path that witnesses $w\in \WR_{\lfloor k/2 \rfloor}[G,L,u]$. Assume to the contrary that $\alpha(u)=\alpha(v)$, then by definition there exists a vertex $w'\in \WR_{\lfloor k/2 \rfloor}[G,L,v]$ with $f(w')=f(w)$ and $sdist(w',v)=sdist(w,u)$. Let $P'(w',v)$ be the path that witnesses $w'\in \WR_{\lfloor k/2 \rfloor}[G,L,v]$, then $|P'(w',v)|=|P(w,u)|$. By taking a path in the union of the path $P(w,v)$ and the path $P'(w',v)$, there exists a path between $w$ and $w'$, which is of length at most $k$, we deduce that one of $w$ and $w'$ is weakly $k$-reachable from the other. As $f(w')=f(w)$, we have $w=w'$ and $k$ is even. As $sdist(w',v)=sdist(w,u)$ and $u,v$ are of distance $k$, we have $|sdist(w',v)|=|sdist(w,u)|=k/2$, which implies that there exists a positive $uv$-path of length $k$, a contradiction.

	\section{Distance 2 in signed planar graphs}\label{sec:planar}
	
	In light of Theorem \ref{stronger}, to prove Theorem \ref{thm:planarexact} it suffices to prove the following result, which is the goal of this section.
	
	\begin{theorem}\label{thm:dcol76}
		For every planar graph $G$  we have $\dcol_4(G)\le 77$.
	\end{theorem}
	
	As mentioned in the introduction, it is possible to get a bound of 76 but, in order to simplify the presentation, we do not do so here. At the end of this section we give a slight hint of how to obtain this better bound.
	 
	Two vertex-disjoint subgraphs $H_1$ and $H_2$ are said to be \emph{adjacent} if there are vertices $v_1\in H_1$ and $v_2\in H_2$ such that $v_{1}v_{2}\in E(G)$.
	A path $P$ in $G$ is called \emph{isometric} if there is no shorter path in $G$ between its endpoints.

	To prove Theorem~\ref{thm:dcol76}, we will create a vertex ordering of $V(G)$ in the following manner: we take an isometric path $P_1$ of $G$ and place its vertices at the start of the ordering of $V(G)$ (ordering the vertices of $P_1$ from one end of the path to the other). Then, we pick a path $P_2$ which is  isometric in $G-P_1$, and put its vertices next in the ordering. We proceed inductively in this way until we have ordered all the vertices. Once we have defined this ordering $L$ of $V(G)$, we will use the following easy lemma to bind $|\DR_{4}[G,L,v]|$ for every $v\in V(G)$.

	\begin{lemma}[Van den Heuvel, Ossona de Mendez, Quiroz, Rabinovich and Siebertz \cite{generalisedcol}] \label{lem:isometricpath}
		Suppose $P$ is an isometric path of $G$ and $v\in V(G)$, then we have $|N^{k}[v]\cap V(P)|\le 2k+1$.
	\end{lemma}

    This motivates the following definition.
	An \emph{isometric-path decomposition} is a set of non-empty paths $\mathcal{P}=\{P_0, \dots, P_s\}$, where $V(G)=\cup_{i=0}^sV(P_i)$, $V(P_i)\cap V(P_j)=\varnothing$ whenever $i\ne j$, and $P_i$ is isometric in $G_i:=G-\cup_{j=0}^{i-1}V(P_j)$ for every $i$. 
	
	We will not be happy with just any isometric path decomposition, but will need one in which every vertex in $G-\cup_{j=0}^{i-1}V(P_j)$ can only be adjacent to a bounded number of paths from $P_0,\dots ,P_{i-1}$. In fact, we will require more properties as follows. (Note that given a graph $G$ and any ordering $L$ of $V(G)$, adding an edge cannot decrease $\DR_4[G,L,v]$. Therefore, to prove Theorem~\ref{thm:dcol76} it is enough to work with maximal planar graphs.)

	\begin{definition}\label{def:reduction}
		A \emph{reduction} of a maximal planar graph $G$ is an isometric-path decomposition $\mathcal{P}=\{P_0, \dots, P_s\}$ of $G$ such that, for a fixed embedding of $G$, we have the following.
		\begin{itemize}
			\item[(i)] The path $P_0$ and $P_1$ are chosen as follows: Let $u_{1}u_{2}u_{3}$ be the vertices in the outerface of $G$, $P_0=u_{1}$ and $P_1=u_2u_3$. Clearly, $G_2$ has only one component, which is in the interior of the cycle $u_1u_2u_3u_1$.
			
			\item[(ii)] For every component $K$ of $G_{i+1}$, the boundary of the face of $G[V(P_0)\cup\cdots\cup V(P_i)]$ which contains $K$ is of the form $C=xP_{h}x'y'P_{j}yx$ for some $h<j\le i$, where $x,x'\in P_h, y,y'\in P_j$ and $xy, x'y'\in E(G)$. In particular, each component of $G_{i+1}$ is adjacent to exactly two paths of $P_0,\dots, P_i$. See Figure~\ref{fig:pathdecomposition}, for an illustration.
			
			\item[(iii)] For all $i\in \{2,\dots, s\}$, $P_i$ is chosen as follows: Let $K_i$ be a non-empty component of $G_{i}$, and let $P_h$ and $P_j$ be the paths adjacent to this component, with $h<j\le i-1$.
			Take two (possibly equal) vertices $w_i, w'_i\in K_i$ such that there are vertices $v_i, v'_i\in P_h$, $z_i, z'_i\in P_j$, with $v_{i}w_{i}z_{i}$ and $v'_{i}w'_{i}z'_{i}$ both bounding faces in $G$ (see again Figure~\ref{fig:pathdecomposition}). Let $P_i$ be a $w_iw_i'$-path which is isometric in $K_i$. 
			If there are multiple choices for $P_i$, then we choose the one that minimizes the number of vertices in the interior of $v_{i}P_{h}v'_{i}w'_{i}P_{i}w_{i}v_{i}$ (the region defined by this cycle which does not contain $P_j$).
		\end{itemize}
	\end{definition}

\begin{figure}[htbp]
		\centering
		\resizebox{.7\textwidth}{!}{
\begin{tikzpicture}[
    line join=round, 
    line cap=round,
    node dot/.style={circle, fill=black, inner sep=1.5pt}
]


\coordinate (Pa0) at (0.5,1.8);
\coordinate (Pa1) at (2,2.5);
\coordinate (Pa2) at (3.5,4.3);
\coordinate (Pa3) at (5,5);
\coordinate (Pa4) at (7.5,5.3);
\coordinate (Pa5) at (9.7,5.3);
\coordinate (Pa6) at (10.8,5.7);

\coordinate (Pb0) at (1,0.3);
\coordinate (Pb1) at (2.5,1);
\coordinate (Pb2) at (4.5,0.5);
\coordinate (Pb3) at (6,0.5);
\coordinate (Pb4) at (7.5,1);
\coordinate (Pb5) at (9.3,2.2);
\coordinate (Pb6) at (10.5,3.6);
\coordinate (Pb7) at (11.5,4);

\coordinate (v1) at (3.5,2);
\coordinate (v2) at (9.5,4.3);
\coordinate (C_left) at (3.5,3);
\coordinate (C_top1) at (4.5,4);
\coordinate (C_top2) at (6,4.3);
\coordinate (C_mid) at (8,4.3);
\coordinate (C_bottom1) at (4.8,1.5);
\coordinate (C_bottom2) at (6.7,1.8);
\coordinate (C_bottom3) at (8.5,3);


\begin{scope}[thin, gray!60]
    \draw (2,2.5) -- (2.5,1);
    \draw (2,2.5) -- (v1);
    \draw (2.5,1) -- (v1);
    \draw (Pa1) -- (C_left);
    \draw (C_left) -- (v1);
    
    \draw (Pa2) -- (C_left);
    \draw (Pa2) -- (v1);
    \draw (Pa2) -- (C_top1);
    \draw (Pa2) -- (C_top2);
    \draw (Pa3) -- (C_top2);
    \draw (C_top2) -- (Pa4);
    \draw (C_top2) -- (C_mid);
    \draw (Pa4) -- (C_mid);
    \draw (Pa5) -- (C_mid);
    
    \draw (Pa5) -- (v2);
    \draw (v2) -- (Pb6);
\end{scope}

\draw[ultra thick] (Pa0) -- (Pa1) -- (Pa2) -- (Pa3) -- (Pa4) -- (Pa5) -- (Pa6);
\draw[ultra thick] (Pb0) -- (Pb1) -- (Pb2) -- (Pb3) -- (Pb4) -- (Pb5) -- (Pb6) -- (Pb7);

\draw[ultra thick] (v1) -- (C_left) -- (C_top1) -- (C_top2) -- (C_mid) -- (v2) -- (C_bottom3);
\draw[ultra thick] (v1) -- (C_bottom1) -- (C_bottom2) -- (C_bottom3);

\draw[thick] (Pa1) -- (Pb1); 
\draw[thick] (Pa5) -- (Pb6); 

\node[node dot, label=below :{\Large $w_i$}] at (v1) {};
\node[node dot, label=below :{\Large $w'_i$}] at (v2) {};
\node[node dot, label=above :{\Large $v'_i$}] at (Pa1) {};
\node[node dot, label=below :{\Large $z'_i$}] at (Pb1) {};
\node[node dot, label=above :{\Large $v_i$}] at (Pa5) {};
\node[node dot, label=below :{\Large $z_i$}] at (Pb6) {};

\node at (4.5,5.5) {\Large $P_h$};
\node at (8,0.5) {\Large $P_j$};
\node at (6,3) {\Large $K_i$};

\end{tikzpicture}

}
		\caption{Component $K_i$ of $G_{i+1}$ is adjacent to $P_h$ and $P_j$. Choice of vertices $w_i, w_i'$ is shown.}
		\label{fig:pathdecomposition}
	\end{figure}

    Van den Heuvel, Ossona de Mendez, Quiroz, Rabinovich and Siebertz \cite{generalisedcol}, used isometric-path decompositions with properties (i), (ii) and (iii), save for the minimality condition of (iii), to prove that for every planar graph $G$ and every integer $k\ge 0$, $\wcol_k(G)\le \binom{k+2}{2}(2k+1)$. Lemma 6.1 of~\cite{generalisedcol} proves the existence of these decompositions for every maximal planar graph. Almulhim and Kierstead \cite{AlmulhimKierstead} introduced the minimality condition of (iii) and the name reduction. This minimality condition can be used to improve the above mentioned bound for $\wcol_k(G)$, for some fixed $k$. In particular, Almulhim and Kierstead used it to improve the upper bound of $\wcol_2(G)$ for $G$ planar, from 30 to 23. If we were to ignore this minimality condition, we would get a bound of 85 in Theorem~\ref{thm:dcol76}, thus not going beyond the bound of 80 given by Alon and Marshall in Theorem~\ref{thm:AM}. Thus this condition is central to our analysis, which inevitably becomes more technical. Reductions as defined here have also been applied in the works of Almulhim \cite{Almulhim} and of Cort\'es, Kumar, Moore, Ossona de Mendez and Quiroz \cite{subchromatic}.

	Before proving Theorem \ref{thm:dcol76} we need some further notation and lemmas. Definition \ref{def:reduction}(ii) tells us that the component $K_i$ is adjacent in $G$ to exactly two paths in $P_0,\dots , P_{i-1}$, but from property~(iii) we further know that each $P_i$ is adjacent to exactly two paths $P_h, P_j$, with $h<j\le i-1$. We call $P_h$ and $P_j$ the \emph{bosses} of $P_i$, with $P_h$ being the \emph{manager} of $P_i$ and $P_j$ being the \emph{foreman} of $P_i$. We denote by $R_{i,1}$ the interior of the cycle $C_{i,1}=v_{i}P_{h}v'_{i}w'_{i}P_{i}w_{i}v_{i}$, that is, the region defined by this cycle which does not contain the foreman $P_j$. Similarly, $R_{i,2}$ is defined to be the interior of the cycle $C_{i,2}=z_{i}P_{j}z'_{i}w'_{i}P_{i}w_{i}z_{i}$.
	For $i\ge 3$, we also set $C_i=v_{i}P_{h}v'_{i}z'_{i}P_{j}z_{i}v_{i}$, that is, $C_i$ is the cycle that is the boundary of the face of $G[V(P_0)\cup\dots \cup V(P_{i-1})]$ which contains $K_i$, and we call the interior of $C_i$ the \emph{region which contains $K_i$}. 

	From now on we assume that $G$ is a maximal planar graph, $\mathcal{P}$ is a fixed reduction of $G$, and~$L$ is an ordering corresponding to $\mathcal{P}$ (as explained at the start of this section). The following observations follow easily from the definition of reduction.

	\begin{observation}\label{obs:intersectManagerForeman}
		Let $v\in P_i$, $i\ge 2$, $P_h$ and $P_j$ be the manager and foreman of $P_i$, respectively. If $P$ is a path between $v$ and a vertex $u\in P_0\cup\cdots\cup P_{i-1}$, then $P$ either intersects $v_{i}P_{h}v'_{i}$ or $z_{i}P_{j}z'_{i}$.
	\end{observation}
    \proof The path $P_i$ lies in the component $K_i$ of $G_i=G-\cup_{j=0}^{i-1}V(P_j)$. But from among $P_0, \dots, P_{i-1}$, the path $P_i$ is adjacent only to $P_h$ and $P_j$. In fact, only $v_{i}P_{h}v'_{i}$ and $z_{i}P_{j}z'_{i}$, by construction. \qed

	\begin{observation}\label{obs:PhisbossofPj}
		Let $P_i\in \mathcal{P}$, $i\ge 2$, $P_h$ and $P_j$ be the manager and foreman of $P_i$, respectively. Then $P_h$ is a boss of $P_j$.
	\end{observation}
    \proof  Since $P_h$ is the manager of $P_i$ we have $h<j$, while property (iii), the bosses of $P_i$ are adjacent through the edges $v_iz_i,v_i'z_i'$. \qed
	
	\begin{observation}\label{obs:subcomponent}
		Let $P_i,P_j\in \mathcal{P}$, with $i\ge 3$. If $P_i$ is a boss of $P_j$, then the component $K_j$ of $G_j$ is a subgraph of the component $K_i$ of $G_i$. 
	\end{observation}
    \proof Being a boss of $P_j$, part of $P_i$ is in the cycle $C_j$, the  boundary of the face of $G[V(P_0)\cup \dots \cup V(P_{j-1})]$ which contains $K_j$. Thus, in $G_i=G-\cup_{\ell\le i-1}V(P_\ell)$ there is an edge joining $V(K_j)$ and $P_i$. But $P_i$ lies in the component $K_i$ of $G_i$. Thus $K_j$ is a (proper) subgraph of $K_i$.\qed
    

	For simplicity, we will write $\dr(v)=\DR_4[G,L,v]$. 
	
	\begin{observation}\label{obs:startinki}
		Let $u\in P_i$ and $v\in P_k$, with $i<k$. Assume that $Q_u=u,u_1,u_2,u_3, v$ is a path witnessing that $u\in \dr(v)$. Then we have $u_3\in K_k$.
	\end{observation}
    \begin{proof}
        By definition of $\dcol_4(G)$, in a path of length 4, the first element, counting from $v$ that can be before $v$ in $L$, is $u_2$. Since $v<_Lu_3$, we cannot have $u_3\in \cup_{j<k}P_j$, and since $vu_3$ is an edge, $u_3$ must be in the same component of $G_k$ as $v$, that is, in $K_k$.
    \end{proof}

	We need the following lemma, which corresponds to Lemma 2.1.5 of \cite{Almulhim} and is a consequence of the minimality condition in point (iii) of the definition of reduction.

	\begin{lemma}[Almulhim \cite{Almulhim}] \label{lem:regionbound}
		Let $P_i\in \mathcal{P}$, $i\ge 2$, and $x,y\in P_i$. Assume that there is an $x y$-path $Q\in R_{i,1}\cup P_i$, such that $Q$ is not contained in $P_i$, then $||xP_iy||\le ||Q||-1$.
	\end{lemma}
	

	\begin{lemma}\label{lem:v in R_i1}
		Given $i\in\{0,\dots ,s\}$, for every $v\in R_{i,1}$, then $|\dr(v)\cap P_i|\le 8$.

	\end{lemma}
	
	\proof Let $x,y\in \dr(v)\cap P_i$ such that $d_{P_i}(x,y)$ is maximum. To prove the lemma, it suffices to prove that $||xP_iy||\le 7$. Let $Q_x$ and $Q_y$ be paths that witnesses $x,y\in \dr(v)$, respectively. In particular, we have that $||Q_x||, ||Q_y||\le 4$, and the walk $W=Q_xvQ_y$ is of length at most 8. As $v$ is in the walk but not on $P_i$, by the minimality of $P_{i}$ given in condition (iii), we have $||xP_iy||\le 7$, as desired.
	\qed

	Theorem \ref{thm:dcol76} is implied by the following result, the proof of which is the aim of the rest of this section. The proof is similar to Theorem 3.0.1 in \cite{Almulhim}, that implies that $\dcol_5(G)\le 95$ for every planar graph $G$.
	
	\begin{theorem}
		We have $|\dr(v)|\le 77$ for every $v\in V(G)$.
	\end{theorem}
	
	\proof 
	
	Towards a contradiction we assume there is a vertex $v$ such that $|\dr(v)|\ge 78$. We let $v\in P_k$ for some $k\ge 0$. If we have $u\in \dr(v)$ and $u\in P_i$, then $i\le k$ by our choice of the vertex ordering. For any $u\in \dr(v)$ we set  $Q_u:=u_{0}\dots u_{q}$ be the path that witnesses $u\in \dr(v)$, where $u_0=u$ and $u_q=v$. 
	
	Recall that $P_i$ is an isometric path in $G_i$. Therefore, by Lemma \ref{lem:isometricpath}, we obtain
    \begin{equation}\label{eq:basic}
    |\dr(v)\cap P_i|\le 9.
    \end{equation}
    Moreover, vertices on $P_k$ that can be reached from $v$ come before $v$ and are at distance at most 4 from $v$ on $P_k$. Hence we have
	\begin{equation}\label{eq:Pk}
		|\dr(v)\cap P_k|\le 5. 
	\end{equation}
	
	Set $\dr^\mathcal{P}(v):=\{P_i \in \mathcal{P}\,|\, \textrm{there exists }u\in  \dr(v)\cap P_i \}$. By \eqref{eq:basic} and \eqref{eq:Pk} we have $|\dr(v)|\le 9\,|\dr^\mathcal{P}(v)|-4$. If $|\dr^\mathcal{P}(v)|\le 8$, which is indeed the case if $k\le 7$, we get $|\dr(v)|\le 68$, which is a contradiction. So we may assume  $k\ge 8$, and in particular $P_k$ has a manager and foreman which we denote by $P_h$ and $P_\ell$, respectively.

	\begin{claim}\label{claim:atmost10}
		We have $|\dr^\mathcal{P}(v)|\le 10$.
	\end{claim}
	
	\proof Let $u\in \dr(v)\setminus\{v\}$. By the definition of $L$ and $\dr(v)$, we have $u<_Lv$ and $u\in V(P_0)\cup\cdots\cup V(P_{k})$. Let $i=\max\{i\in \{0,\dots, q-1\}: u_i<_L v\}$ (recall $Q_u=u_{0}\dots u_{q}$, with $u_0=u$). By the definition of $\dr(v)$, we have $i\le 2$. We will take cases on the value of $i$ to see on which path can $u$ lie.
	
	Case 1: $i=0$.
	Assume that $u\notin P_k$, then by Observation \ref{obs:intersectManagerForeman}, $Q_{u}\cap (P_h\cup P_\ell)\ne \varnothing$. Since $u$ is the only vertex in $Q_u$ which is smaller than $v$ with respect to $L$, we have $u\in P_h\cup P_\ell$. Together with the fact that $u$ might be in $P_k$, $u\in P_k\cup P_\ell\cup P_h$.
	
	Case 2: $i=1$.
	The argument of Case 1 gives us that $u_1\in P_k\cup P_\ell\cup P_h$. By the definition of $\dr(v)$, $u<_Lu_1$. If $u_1\in P_k$, then by Observation \ref{obs:intersectManagerForeman}, $u\in P_k\cup P_\ell\cup P_h$. First assume we have $u_1\in P_h$. Since $u<_Lu_1$, we must have a vertex which is smaller than $u_1$, and thus, since $|V(P_0)|=1$, we must have $h\ge 1$. Let $P_g$ and $P_f$ be the manager and foreman of $P_h$, respectively (where we take $P_f=P_g=P_0$, if $h=1$). Then, by Observation~\ref{obs:intersectManagerForeman}, $u\in P_g\cup P_f\cup P_h$. Finally assume $u_1\in P_\ell$ with $\ell\ge 1$. Since $P_h$ and $P_\ell$ are the bosses of $P_k$ with $h<\ell$, by Observation~\ref{obs:PhisbossofPj}, $P_h$ is a boss of $P_\ell$. Let $P_m$ be the other boss of $P_\ell$. Again Observation~\ref{obs:intersectManagerForeman} gives us that $u\in P_\ell\cup P_m\cup P_h$. Thus when $i=1$, we have $u\in P_j$ for some $j\in A:=\{k,\ell,m,h,f,g\}$.
	
	Case 3: $i=2$.
	The argument of Case 2 gives us that $u_1\in P_j$ for some $j\in A$. If $u_1\in P_k\cup P_\ell\cup P_h$, then the argument of Case 2 gives us that $u\in P_j$ for some $j\in A$. So we can take $u_1\in P_j$, with $j\in\{m,f,g\}$. First assume that $u_1\in P_g$ for $g\ge 1$. Let $P_b$ and $P_c$ be the manager and foreman of $P_g$, respectively. By Observation~\ref{obs:intersectManagerForeman}, $u\in P_g\cup P_b\cup P_c$. Now assume that $u_1\in P_f$ with $f\ge 1$. Recall that $P_g$ and $P_f$ are the manager and foreman of $P_h$, so by Observation~\ref{obs:PhisbossofPj}, $P_g$ is a boss of $P_f$. Let $P_d$ be the other boss of $P_f$. By Observation~\ref{obs:intersectManagerForeman}, $u\in P_f\cup P_d\cup P_g$.
	
	Assume that $u_1\in P_m$. Recall that $P_m$ and $P_h$ are the bosses of $P_\ell$. By Observation~\ref{obs:PhisbossofPj}, either $P_m$ is a boss of $P_h$ or vice versa. If $P_m$ is a boss of $P_h$, then $m\in \{f, g\}$ and we have already considered this possibility. So assume that $P_h$ is a boss of $P_m$ and thus $m\notin\{f, g\}$. Let $P_n$ be the other boss of $P_m$. By Observation~\ref{obs:intersectManagerForeman}, $u\in P_m\cup P_n\cup P_h$. So when $i=2$, (and also when $i\le 1$) $u\in P_j$ for some $j\in B:=\{k,\ell, m,n,h,f,d,g,c,b\}$; see the example in Figure~\ref{fig:F1}. \hfill $\Diamond$
	

	\begin{figure}[htbp]
		\centering
		\resizebox{.7\textwidth}{!}{
			\begin{tikzpicture}[scale=.4]
				
				
				\draw(-1, 4) node[circle, draw=black!80, inner sep=0mm, minimum size=.8mm, fill=black] (k1){};
				
				\foreach \i in {5,4,3,2}
				{
					\draw[yshift=-21*(\i-1)] (-1, 4) node[circle, draw=black!80, inner sep=0mm, minimum size=.8mm, fill=black] (k\i){};
				}
				
				\draw (-.4, 2.6) node[label= $P_k$]  (Pk){};
				\draw (-.5, 1.8) node[label= $v$]  (v){};
				
				\draw  [line width=0.3mm, black] (k1) -- (k5);
				
				
				\draw(1, 6) node[circle, draw=black!80, inner sep=0mm, minimum size=.8mm, fill=black] (l1){};
				
				\foreach \i in {6,5,4,3,2}
				{
					\draw[yshift=-40*(\i-1)] (1, 6) node[circle, draw=black!80, inner sep=0mm, minimum size=.8mm, fill=black] (l\i){};
				}
				
				\draw (1.6, 1.8) node[label= $P_\ell$]  (Pl){};
				
				\draw  [line width=0.3mm, black] (l1) -- (l6);
				
				
				\draw(-3, 12) node[circle, draw=black!80, inner sep=0mm, minimum size=.8mm, fill=black] (h1){};
				
				\foreach \i in {4,3,2}
				{
					\draw[yshift=-50*(\i-1)] (-3,12) node[circle, draw=black!80, inner sep=0mm, minimum size=.8mm, fill=black] (h\i){};
				}

				\draw(-3, -7) node[circle, draw=black!80, inner sep=0mm, minimum size=.8mm, fill=black] (h11){};
				\foreach \i in {10,9,8}
				{
					\draw[yshift=-50*(\i-11)] (-3,-7) node[circle, draw=black!80, inner sep=0mm, minimum size=.8mm, fill=black] (h\i){};
				}

				\node[circle, draw=black!80, inner sep=0mm, minimum size=.8mm, fill=black, shift={(-1.6,0)}] (h5) at (l2) {};
				\node[circle, draw=black!80, inner sep=0mm, minimum size=.8mm, fill=black, shift={(-1.6,0)}] (h7) at (l5) {};
				\draw(-3, 2.6) node[circle, draw=black!80, inner sep=0mm, minimum size=.8mm, fill=black] (h6){};

				\draw (-2.2, 1.8) node[label= $P_h$]  (Ph){};
				
				\draw  [line width=0.3mm, black] (h1) -- (h11);
				\draw  [line width=0.3mm, green] (h5) -- (l2) (h5) -- (k1) (l2) -- (k1);
				\draw  [line width=0.3mm, green] (h7) -- (l5) (h7) -- (k5) (l5) -- (k5);
				
				\node[circle, draw=black!80, inner sep=0mm, minimum size=.8mm, fill=black, label={above: $w_h$}] (h0) at (h1) {};
				
				\node[circle, draw=black!80, inner sep=0mm, minimum size=.8mm, fill=black, label={below, yshift=.1cm: $w'_h$}] (h12) at (h11) {};

				\node[circle, draw=black!80, inner sep=0mm, minimum size=.8mm, fill=black, shift={(2.8,0)}] (m2) at (h4) {};
				\node[circle, draw=black!80, inner sep=0mm, minimum size=.8mm, fill=black, shift={(0,.4)}] (m1) at (m2) {};

				\foreach \i in {3,4,5,6}
				{
					\node[circle, draw=black!80, inner sep=0mm, minimum size=.8mm, fill=black, shift={(0,.67*(\i-7))}] (m\i) at (m2) {};
				}
				
				\node[circle, draw=black!80, inner sep=0mm, minimum size=.8mm, fill=black, shift={(2.8,0)}] (m7) at (h8) {};
				\node[circle, draw=black!80, inner sep=0mm, minimum size=.8mm, fill=black, shift={(0,-.4)}] (m8) at (m7) {};
				
				\draw (4.7, 1.8) node[label= $P_m$]  (Pm){};
				
				\draw  [line width=0.3mm, black] (m1) -- (m8);
				\draw  [line width=0.3mm, green] (m2) -- (h4) (h4) -- (l1) (m2) -- (l1);
				\draw  [line width=0.3mm, green] (m7) -- (h8) (h8) -- (l6) (l6) -- (m7);
				
				
				\node[circle, draw=black!80, inner sep=0mm, minimum size=.8mm, fill=black, shift={(4,0)}] (n1) at (h2) {};
				\node[circle, draw=black!80, inner sep=0mm, minimum size=.8mm, fill=black, shift={(4,0)}] (n2) at (h3) {};
				\node[circle, draw=black!80, inner sep=0mm, minimum size=.8mm, fill=black, shift={(4,0)}] (n7) at (h9) {};
				\node[circle, draw=black!80, inner sep=0mm, minimum size=.8mm, fill=black, shift={(4,0)}] (n8) at (h10) {};
				
				\foreach \i in {3,4,5,6}
				{
					\node[circle, draw=black!80, inner sep=0mm, minimum size=.8mm, fill=black, shift={(0,-.95*(\i-2))}] (n\i) at (n2) {};
				}
				
				\draw  [line width=0.3mm, black] (n1) -- (n8);
				
				\draw  [line width=0.3mm, green] (n1) -- (h2) (n2) -- (h3) (n2) -- (m1) (m1) -- (h3);
				\draw  [line width=0.3mm, green] (n8) -- (h10) (n7) -- (h9) (n7) -- (m8) (m8) -- (h9);
				
				\draw (7.7, 1.8) node[label= $P_n$]  (Pn){};
				
				
				\node[circle, draw=black!80, inner sep=0mm, minimum size=.8mm, fill=black, shift={(-1.5,.5)}] (f2) at (h1) {};
				\node[circle, draw=black!80, inner sep=0mm, minimum size=.8mm, fill=black, shift={(-1.5,-.5)}] (f7) at (h11) {};
				\node[circle, draw=black!80, inner sep=0mm, minimum size=.8mm, fill=black, shift={(0,.55)}] (f1) at (f2) {};
				\node[circle, draw=black!80, inner sep=0mm, minimum size=.8mm, fill=black, shift={(0,-.55)}] (f8) at (f7) {};
				
				\draw  [line width=0.3mm, black] (f1) -- (f8);
				
				\foreach \i in {3,4,5,6}
				{
					\node[circle, draw=black!80, inner sep=0mm, minimum size=.8mm, fill=black, shift={(0,-1.7*(\i-2))}] (f\i) at (f2) {};
				}
				
				\draw (-6, 1.8) node[label= $P_f$]  (Pf){};
				
				\node[circle, draw=black!80, inner sep=0mm, minimum size=.8mm, fill=black, label={left: $z_h$}] (h9) at (f2) {};
				
				\node[circle, draw=black!80, inner sep=0mm, minimum size=.8mm, fill=black, label={left: $z'_h$}] (h10) at (f7) {};
				
				
				\node[circle, draw=black!80, inner sep=0mm, minimum size=.8mm, fill=black, shift={(7,0)}] (g3) at (f2) {};
				\node[circle, draw=black!80, inner sep=0mm, minimum size=.8mm, fill=black, shift={(7,0)}] (g6) at (f7) {};
				
				\foreach \i in {2}
				{
					\node[circle, draw=black!80, inner sep=0mm, minimum size=.8mm, fill=black, shift={(0,1*(3-\i))}] (g\i) at (g3) {};
				}
				
				\foreach \i in {7}
				{
					\node[circle, draw=black!80, inner sep=0mm, minimum size=.8mm, fill=black, shift={(0,1*(6-\i))}] (g\i) at (g6) {};
				}
				
				\node[circle, draw=black!80, inner sep=0mm, minimum size=.8mm, fill=black, shift={(0,.6)}] (g1) at (g2) {};
				\node[circle, draw=black!80, inner sep=0mm, minimum size=.8mm, fill=black, shift={(0,-.6)}] (g8) at (g7) {};

				\draw  [line width=0.3mm, black] (g1) -- (g8);
				
				\foreach \i in {4,5}
				{
					\node[circle, draw=black!80, inner sep=0mm, minimum size=.8mm, fill=black, shift={(0,3*(3-\i))}] (g\i) at (g3) {};
				}
				
				\draw (11.4, 1.8) node[label= $P_g$]  (Pg){};
				
				\draw  [line width=0.3mm, green] (g3) -- (h1) (h1) -- (f2) (f2) -- (g3);
				\draw  [line width=0.3mm, green] (g6) -- (h11) (h11) -- (f7) (f7) -- (g6);
				
				\node[circle, draw=black!80, inner sep=0mm, minimum size=.8mm, fill=black, label={right: $v_h$}] (g9) at (g3) {};
				
				\node[circle, draw=black!80, inner sep=0mm, minimum size=.8mm, fill=black, label={right: $v'_h$}] (g10) at (g6) {};
				
				
				\node[circle, draw=black!80, inner sep=0mm, minimum size=.8mm, fill=black, shift={(-8.5,0)}] (d1) at (g2) {};
				\node[circle, draw=black!80, inner sep=0mm, minimum size=.8mm, fill=black, shift={(-8.5,0)}] (d5) at (g7) {};
				
				\draw  [line width=0.3mm, black] (d1) -- (d5);
				
				\foreach \i in {2,3,4}
				{
					\node[circle, draw=black!80, inner sep=0mm, minimum size=.8mm, fill=black, shift={(0,2.5*(1-\i))}] (d\i) at (d1) {};
				}
				
				\draw (-9.8, 1.8) node[label= $P_d$]  (Pd){};
				
				\draw  [line width=0.3mm, green] (d1) -- (f1) (f1) -- (g2) (g2) -- (d1);
				\draw  [line width=0.3mm, green] (d5) -- (g7) (g7) -- (f8) (f8) -- (d5);
				
				
				\node[circle, draw=black!80, inner sep=0mm, minimum size=.8mm, fill=black, shift={(1.9,.6)}] (c1) at (g1) {};
				\node[circle, draw=black!80, inner sep=0mm, minimum size=.8mm, fill=black, shift={(1.9,-.6)}] (c4) at (g8) {};
				
				\draw  [line width=0.3mm, black] (c1) -- (c4);
				
				\foreach \i in {2,3}
				{
					\node[circle, draw=black!80, inner sep=0mm, minimum size=.8mm, fill=black, shift={(0,4.6*(1-\i))}] (c\i) at (c1) {};
				}
				
				\draw (16.4, 1.8) node[label= $P_c$]  (Pc){};
				
				
				\node[circle, draw=black!80, inner sep=0mm, minimum size=.8mm, fill=black, shift={(-12.3,0)}] (b1) at (c1) {};
				\node[circle, draw=black!80, inner sep=0mm, minimum size=.8mm, fill=black, shift={(-12.3,0)}] (b4) at (c4) {};
				
				\draw  [line width=0.3mm, black] (b1) -- (b4);
				
				\foreach \i in {2,3}
				{
					\node[circle, draw=black!80, inner sep=0mm, minimum size=.8mm, fill=black, shift={(0,4.6*(1-\i))}] (b\i) at (b1) {};
				}
				
				\draw (-14.4, 1.8) node[label=$P_b$]  (Pb){};
				
				\draw  [line width=0.3mm, green] (b1) -- (c1) (b1) -- (g1) (g1) -- (c1);
				\draw  [line width=0.3mm, green] (b4) -- (c4) (b4) -- (g8) (g8) -- (c4);

			\end{tikzpicture}
		}
		\caption{$v$ belongs to $R_{h,1}$ which is the interior of $C_{h,1}=v_{h}P_{g}v'_{h}w'_{h}P_{h}w_{h}v_{h}$.}
		\label{fig:F1}
	\end{figure}

	From now on, we use the notation introduced in the proof of the previous claim by which we have that if $P_i\in \dr^\mathcal{P}(v)$, then $i\in\{k,\ell,m,n,h,f,d,g,c,b\}$. Moreover it gives us the relationships summarized in Table \ref{table:bad}.

    \begin{table}[H]
	\begin{center}
		\begin{tabular}{|c|c|c|}
			\hline
			Path&	Bosses & Manager (if set)\\
			\hline
			$P_k$& $P_\ell, P_h$& $P_h$\\
			
			$P_\ell$& $P_m, P_h$& \\
			
			$P_h$& $P_f, P_g$& $P_g$\\
			
			$P_g$& $P_b, P_c$& $P_b$\\
			
			$P_f$& $P_g, P_d$&\\
			
			$P_m$& $P_h, P_n$&\\
			\hline
		\end{tabular}	
		\caption{Relationships between the paths in $\dr^\mathcal{P}(v).$}
        \label{table:bad}
	\end{center}
    \end{table}
    
	So far we have shown that $|\dr^\mathcal{P}(v)|\in \{9,10\}$, our next task is to show that this quantity is exactly 10.

	\begin{claim}\label{claim:allmostordered}
		We have $h<m$.
	\end{claim}
	
	\proof 	Towards a contradiction take $m<h$. Recall that $P_h$ and $P_m$ are the bosses of $P_\ell$. Then by Observation \ref{obs:PhisbossofPj}, $P_m$ is a boss of $P_h$. Since Table~\ref{table:bad} tells us that $P_f$ and $P_g$ are the bosses of $P_h$, we have that $m\in \{f,g\}$. The last two rows of Table \ref{table:bad} tell us that  if $m=f$, then we have $\{h,n\}=\{g,d\}$. Similarly, if $m=g$, then we have $\{h,n\}=\{b,c\}$. Since we have only 10 paths that can belong to $\dr^\mathcal{P}(v)$, any pair of identifications contradicts the fact that $|\dr^\mathcal{P}(v)|\in \{9,10\}$. Therefore, we must have $h<m$.
	\hfill $\Diamond$

	In particular, we now have that the manager and foreman of $P_\ell$ are $P_h$ and $P_m$, respectively.
	
	\begin{claim}\label{claim:ell}
		We have $|\dr(v)\cap P_\ell|\le 8$. 
	\end{claim}
	
	\proof 
	As $P_h$ and $P_m$ are the bosses of $P_\ell$ and, by Claim~\ref{claim:allmostordered}, $h<m$, then $P_h$ is the manager of $P_\ell$. But since $P_h$ is also the manager of $P_k$, we have $v\in R_{\ell,1}$. We may now apply Lemma~\ref{lem:v in R_i1}. 
	\hfill $\Diamond$
	
	\begin{claim}\label{claim:drp=10}
		We have $|\dr^\mathcal{P}(v)|=10$.
	\end{claim}
	
	\proof
	Observe that $P_k, P_\ell\in \dr^\mathcal{P}(v)$. If $|\dr^\mathcal{P}|=9$, then by \eqref{eq:basic}, \eqref{eq:Pk} and Claim \ref{claim:ell} we have  $|\dr^\mathcal{P}|\le 76$, a contradiction. \phantom{word}\hfill $\Diamond$

	\begin{claim}\label{claim:finishordering}
		We have $g<d$ and $h<n$.
	\end{claim}
	
	\proof Suppose that $d<g$. Since $P_g$ and $P_d$ are the bosses of $P_f$, the manager and foreman of $P_f$ are $P_d$ and $P_g$, respectively. By Observation~\ref{obs:PhisbossofPj}, $P_d$ is a boss of $P_g$, that is, $d\in \{b, c\}$. This contradicts Claim \ref{claim:drp=10}.
	
	Proof of  $h<n$, follows similarly. 
	\hfill $\Diamond$

	With Claims~\ref{claim:allmostordered} and \ref{claim:finishordering}, we can fill the missing information in the Table~\ref{table:bad} as follows:
	
	\begin{table}[H]
	\begin{center}
		\begin{tabular}{|c|c|c|}
			\hline
			Path&	Bosses & Manager\\
			\hline
			
			$P_\ell$& $P_m, P_h$&$P_h$ \\
			
			$P_f$& $P_g, P_d$&$P_g$\\
			
			$P_m$& $P_h, P_n$&$P_h$\\
			\hline
			
		\end{tabular}	
		\caption{Supplement of Table~\ref{table:bad}}
        \label{table:excellent}
	\end{center}
    \end{table}

	We have that $$\dr^\mathcal{P}(v) = \{ k,\ell , m,n,h,f,d,g,c,b\},$$ and by the information in Tables \ref{table:bad} and \ref{table:excellent} we can conclude that we have
	$$b<c<g<d<f<h<n<m<\ell<k.$$ 
(Note that Figure~\ref{fig:F1} is compatible with this ordering of the paths.)     In particular, now that we know that $P_g$ is the manager of $P_f$, we can argue in a similar way as in Claim~\ref{claim:ell} to prove  the following.
	
	\begin{claim}\label{claim:f}
		We have $|\dr(v)\cap P_f|\le 8$.
	\end{claim}
		

We now prove bounds for $P_m$ and $P_n$.
        
		\begin{claim}\label{claim:pathm}
			We have $|\dr(v)\cap P_m| \le 7$.
		\end{claim}
        \proof Given $x,y \in \dr(v)\cap P_m$, it suffices to prove that $||xP_my||\le 6$. Let $Q_x$ and $Q_y$ be paths that witness $x,y\in \dr(v)$, respectively. If one of $Q_x, Q_y$ is of length less than 4, then $Q_x\cup Q_y$ is a walk of length at most 7. We have $v\in Q_x\cup Q_y$ and $v\in R_{m,1}$, but also that $v$ is not in $P_m$. Therefore, by the minimality of $P_m$ given by condition (iii), we have $||xP_my||\le 6$. Thus we may assume $||Q_x||=||Q_y||=4$. 
        
        Since $P_h$ and $P_m$ are the manager and foreman of $P_\ell$, we have $h<m$ and therefore $Q_x\cap P_h$, $Q_y\cap P_h=\varnothing$. Furthermore, $v\in R_{\ell,1}$, so $Q_x\cap P_\ell$, $Q_y\cap P_\ell\neq \varnothing$. Starting from $v$, let $x'$ and $y'$ be the first vertex in $Q_x$ and $Q_y$, respectively, contained in $P_\ell$. Then the walk $W:=x'Q_xvQ_yy'$ is in $R_{\ell,1}$. By Lemma~\ref{lem:regionbound}, $||x'P_\ell y'||\le ||W||-1$. As $P_h$ is the manager of $P_m$, and $P_m$ is the foreman of $P_\ell$, the walk $W':=xQ_xx'P_\ell y'Q_yy$ is in $R_{m,1}$. Again by Lemma~\ref{lem:regionbound}, 
	\begin{align*}
	||xP_my||&\le ||W'||-1\\
	&=||xQ_xx'||+||x'P_\ell y'||+||y'Q_yy||-1\\
	&\le ||xQ_xx'||+||W||-1+||y'Q_yy||-1\\
	&= ||xQ_xvQ_yy||-2\\
	&\le 6, 
	\end{align*}
	as desired.\hfill $\Diamond$

		By a similar argument we can obtain the following.
		
		\begin{claim}\label{claim:pathn}
			We have $ |\dr(v)\cap P_n|\le 7$.
		\end{claim}

		Our last goal is to prove that we have 
		\begin{equation}\label{eq:dh}
			|\dr(v)\cap(P_d\cup P_h)|\le 15.
		\end{equation}
		For this we need the five claims below.
		
		\begin{claim}\label{claim:towardsd}
			Let $u\in \dr(v)\cap P_d$ and assume that $Q_u:=u_{0}u_{1}u_{2}\dots u_{q}$ is a path witnesses that $u\in \dr(v)$ where $u_0=u$ and $u_q=v$. Then we have $u_1\in P_f$ and $u_2\in v_kP_hv'_k$.
		\end{claim}
		
		\proof Let $i=\max\{j\in \{0,\dots,q-1\}\mid u_j<_L v \}$. From the proof of Claim~\ref{claim:atmost10}, we have $u\in \{P_h, P_k, P_\ell\}$ if $i=0$, and that $u\in P_j$ for some $j\in \{k, \ell,m,h,f,g\}$ if $i=1$. Thus we cannot have $i\le 1$. Moreover, when $i=2$, we cannot have $u_1\in \{P_h, P_k, P_\ell\}$, as otherwise we would again have $u\in P_j$ for some $j\in \{k, \ell,m,h,f,g\}$. Thus we must have $u_1\in \{P_f, P_g, P_m\}$. The proof of Claim \ref{claim:atmost10} also gives that we only have $u_0\in P_d$ when we have $u_1\in P_f$ (or when $u\in P_m$ and $P_m=P_f$, but we already have $f<m$). 
        
        From Observation \ref{obs:intersectManagerForeman}, the path $u_1Q_uv$ intersects with either $v_kP_hv'_k$ or $P_\ell$. So either $u_2\in v_kP_hv'_k$ or $u_2\in P_\ell$. Since $P_f$ and $P_\ell$ are not adjacent, we obtain $u_2\in v_kP_hv'_k$.
		\hfill $\Diamond$
		
		\begin{claim}\label{claim:pathd}
			$|\dr(v)\cap P_d|\le 8$.
		\end{claim}
		
		\proof
		Consider vertices $x,y\in \dr(v)\cap P_d$, and paths $Q_x$ and $Q_y$ witnessing that $x$ and $y$, respectively, belong to $\dr(v).$ Note that we have  $3\le ||Q_x||, ||Q_y||\le 4$. From Claim~\ref{claim:towardsd}, if $Q_x=x,x',x'',\dots, v$ and $Q_y=y,y',y'',\dots, v$, then $x',y'\in P_f$ and $x'',y''\in v_kP_hv_k'$. We want to use Lemma~\ref{lem:regionbound} and for this we show that the walk $W:= x'Q_xvQ_yy'$ lies in $P_f\cup R_{f,1}$. Since $P_g$ and $P_f$ are the bosses of $P_h$ we have that $K_h\subseteq R_{f,1}$. From Observation~\ref{obs:startinki}, we have $x''Q_xv-x'',y''Q_yv-y''\subseteq K_k$, and by Observation~\ref{obs:subcomponent}, $K_k\subset K_h$. Thus we have $x''Q_xvQ_y''\subseteq K_h$. We mentioned that $K_h\subseteq R_{f,1}$,
		and since we have $x',y'\in P_f$, we obtain $W\subseteq P_f\cup R_{f,1}$. Therefore we can use Lemma~\ref{lem:regionbound} to get $||x'P_fy'||\le ||W||-1$.  
		
		Since $x,y\in P_d$ and $x'P_fy'\subseteq K_f\subset K_d$, we have that the walk $xx'P_fy'y$ is contained in $K_d$. As $xP_dy$ is isometric in $K_d$, we have
		\begin{align*}
			||xP_dy||&\le ||xx'P_fy'y||\\
			&=||xx'||+||x'P_fy'||+||y'y||\\
			&\le ||xx'||+||W||+||y'y||-1\\
			&\le ||xQ_xvQ_yy||-1\\
			&\le 7, 
		\end{align*}
		as desired. \hfill $\Diamond$
		
		
		\begin{claim}\label{claim:pathdrh1}
			If $v\in R_{h,1}$, then we have $|\dr(v)\cap P_d|\le 7$.
		\end{claim}
		
		\proof 
		
		Consider vertices $x,y\in \dr(v)\cap P_d$, and paths $Q_x$ and $Q_y$ witnessing that $x$ and $y$, respectively, belong to $\dr(v).$ We have  $||Q_x||, ||Q_y||\le 4$ and $Q_x,Q_y\subseteq K_d$. By Claim~\ref{claim:towardsd} , if $Q_x=x,x',x'',\dots, v$ and $Q_y=y,y',y'',\dots, v$, then $x',y'\in P_f$ and $x'',y''\in v_kP_hv_k'$. By Observation~\ref{obs:startinki} the walk $W:=x''Q_xvQ_yy''$ satisfies $W-\{x'',y''\}\subseteq K_k$. By Observation~\ref{obs:subcomponent}, we have $K_k\subset K_h$, and thus either $K_k\subseteq R_{h,1}$ or $K_k\subseteq R_{h,2}$. But by hypothesis we have $v\in R_{h,1}$, and thus $K_k\subseteq R_{h,1}$. We obtain $W\subseteq P_h\cup R_{h,1}$ and by Lemma~\ref{lem:regionbound}, we obtain $||x''P_hy''||\le ||W||-1$.

		Recall that $C_{f,1}=v_fP_gv'_fw'_fP_fw_fv_f$ is the cycle that forms the frontier of $R_{f,1}$. As $P_g$ and $P_f$ are the bosses of $P_h$, and $P_g$ is the manager of $P_f$, $P_h$ is in the interior of $C_{f,1}$. Thus the walk $W':=x'x''P_hy''y'$ lies in $P_f\cup R_{f,1}$, and by Lemma~\ref{lem:regionbound}, we have
		\begin{align*}
			||x'P_fy'||&\le ||W'||-1\\
			&=||x'x''||+||x''P_hy''||+||y''y'||-1\\
			&\le ||x'x''||+||W||+||y''y'||-2\\
			&= ||x'Q_xvQ_yy'||-2.
		\end{align*}
		
		By Observation~\ref{obs:subcomponent}, we have $x'P_fy'\subseteq K_f\subset K_d$ and, since $x,y\in P_d$, the walk $W'':=xx'P_fy'y$ is in $K_d$. Since $xP_dy$ is isometric in $K_d$, we have
		\begin{align*}
			||xP_dy||&\le ||W''||\\
			&=||xx'||+||x'P_fy'||+||y'y||\\
			&\le ||xx'||+||x'Q_xvQ_yy'||+||y'y||-2\\
			&=||xQ_xvQ_yy||-2\\
			&\le 6,
		\end{align*}
		as desired. \hfill $\Diamond$
		
		\begin{figure}[htbp]
			\centering
			\resizebox{.6\textwidth}{!}{
				\begin{tikzpicture}[scale=.4]
					
					
					\draw(1, 4) node[circle, draw=red!80, inner sep=0mm, minimum size=.8mm, fill=red,label={above: $w_k$}] (h1){};
					
					\node[circle, draw=red!80, inner sep=0mm, minimum size=.8mm, fill=red, shift={(0,-.5)}] (h2) at (h1) {};

					\node[circle, draw=red!80, inner sep=0mm, minimum size=.4mm, fill=red, shift={(0,-.5)}] (h3) at (h2) {};
					
					\foreach \i in {4,5}
					{
						\node[circle, draw=red!80, inner sep=0mm, minimum size=.4mm, fill=red, shift={(0,.25*(3-\i))}] (h\i) at (h3) {};
					}
					\node[circle, draw=red!80, inner sep=0mm, minimum size=.8mm, fill=red, shift={(0,-.5)}] (h6) at (h5) {};

					\node[circle, draw=red!80, inner sep=0mm, minimum size=.8mm, fill=red, shift={(0,-.5)}, label={below,yshift=.08cm: $w'_k$}] (h7) at (h6) {};
					
					\node[circle, draw=red!80, inner sep=0mm, minimum size=.4mm, fill=red, label={right: $P_k$}] (h8) at (h4) {};

					\draw  [line width=0.3mm, red] (h1) -- (h2) (h7) -- (h6);
					
					
					\node[circle, draw=black!80, inner sep=0mm, minimum size=.8mm, fill=black, shift={(-2,0)}] (f5) at (h2) {};
					\node[circle, draw=black!80, inner sep=0mm, minimum size=.8mm, fill=black, shift={(-2,0)}] (f9) at (h6) {};
					\node[circle, draw=black!80, inner sep=0mm, minimum size=.8mm, fill=black, shift={(0,1)},label={left: $z_k$}] (f4) at (f5) {};
					\node[circle, draw=black!80, inner sep=0mm, minimum size=.8mm, fill=black, shift={(0,-1)}, label={left: $z'_k$}] (f10) at (f9) {};
					
					\node[circle, draw=black!80, inner sep=0mm, minimum size=.4mm, fill=black, shift={(0,-.5)}] (f6) at (f5) {};
					
					\foreach \i in {7,8}
					{
						\node[circle, draw=black!80, inner sep=0mm, minimum size=.4mm, fill=black, shift={(0,.25*(6-\i))}] (f\i) at (f6) {};
					}
					
					\node[circle, draw=black!80, inner sep=0mm, minimum size=.4mm, fill=black, label={right: $P_\ell$}] (h22) at (f7) {};
					
					\node[circle, draw=black!80, inner sep=0mm, minimum size=.4mm, fill=black, shift={(0,.5)}] (f3) at (f4) {};
					\foreach \i in {1,2}
					{
						\node[circle, draw=black!80, inner sep=0mm, minimum size=.4mm, fill=black, shift={(0,.25*(3-\i))}] (f\i) at (f3) {};
					}
					
					\node[circle, draw=black!80, inner sep=0mm, minimum size=.4mm, fill=black, shift={(0,-.5)}] (f11) at (f10) {};
					\foreach \i in {12,13}
					{
						\node[circle, draw=black!80, inner sep=0mm, minimum size=.4mm, fill=black, shift={(0,.25*(11-\i))}] (f\i) at (f11) {};
					}
					
					\draw  [line width=0.3mm, black] (f4) -- (f5) (f9) -- (f10);
					
					
					\node[circle, draw=blue!80, inner sep=0mm, minimum size=.8mm, fill=blue, shift={(2,0)}] (g8) at (h2) {};
					\node[circle, draw=blue!80, inner sep=0mm, minimum size=.8mm, fill=blue, shift={(2,0)}] (g12) at (h6) {};
					\node[circle, draw=blue!80, inner sep=0mm, minimum size=.8mm, fill=blue, shift={(2,1)}, label={right: $v_k=x''$}] (g7) at (h2) {};
					\node[circle, draw=blue!80, inner sep=0mm, minimum size=.8mm, fill=blue, shift={(2,-1)}, label={right: $v'_k=y''$}] (g13) at (h6) {};
					
					\draw  [line width=0.3mm, blue] (g7) -- (g8) (g12) -- (g13);
					
					\node[circle, draw=blue!80, inner sep=0mm, minimum size=.4mm, fill=blue, shift={(0,-.5)}] (g9) at (g8) {};
					\foreach \i in {10,11}
					{
						\node[circle, draw=blue!80, inner sep=0mm, minimum size=.4mm, fill=blue, shift={(0,.25*(9-\i))}] (g\i) at (g9) {};
					}
					
					\node[circle, draw=blue!80, inner sep=0mm, minimum size=.4mm, fill=blue, label={right: $P_h$}] (g22) at (g10) {};

					\draw  [line width=0.3mm, green] (f4) -- (h1) (f4) -- (g7) (f10) -- (h7) (f10) -- (g13);
					\draw  [line width=0.3mm, red] (g7) -- (h1) (h7) -- (g13);
					
					\node[circle, draw=red!80, inner sep=0mm, minimum size=.8mm, fill=red, shift={(0,.5)}] (g6) at (g7) {};
					\node[circle, draw=red!80, inner sep=0mm, minimum size=.8mm, fill=red, shift={(0,-.5)}] (g14) at (g13) {};

					\draw  [line width=0.3mm, red] (g6) -- (g7) (g13) -- (g14);
					
					\node[circle, draw=red!80, inner sep=0mm, minimum size=.4mm, fill=red, shift={(0,.5)}] (g5) at (g6) {};
					\foreach \i in {3,4}
					{
						\node[circle, draw=red!80, inner sep=0mm, minimum size=.4mm, fill=red, shift={(0,.25*(5-\i))}] (g\i) at (g5) {};
					}
					
					\foreach \i in {1,2}
					{
						\node[circle, draw=red!80, inner sep=0mm, minimum size=.8mm, fill=red, shift={(0,.5*(3-\i))}] (g\i) at (g3) {};
					}
					\node[circle, draw=red!80, inner sep=0mm, minimum size=.8mm, fill=red, shift={(0,.5)},label={above: $w_h$}] (g1) at (g2) {};

					\node[circle, draw=red!80, inner sep=0mm, minimum size=.4mm, fill=red, shift={(0,-.5)}] (g15) at (g14) {};
					\foreach \i in {16,17}
					{
						\node[circle, draw=red!80, inner sep=0mm, minimum size=.4mm, fill=red, shift={(0,.25*(15-\i))}] (g\i) at (g15) {};
					}
					
					\foreach \i in {18,19}
					{
						\node[circle, draw=red!80, inner sep=0mm, minimum size=.8mm, fill=red, shift={(0,.5*(17-\i))}] (g\i) at (g17) {};
					}
					\node[circle, draw=red!80, inner sep=0mm, minimum size=.8mm, fill=red, shift={(0,-.5)},label={below,yshift=.08cm: $w'_h$}] (g19) at (g18) {};
					
					\draw  [line width=0.3mm, red] (g1) -- (g2) (g18) -- (g19);
					
					
					\node[circle, draw=red!80, inner sep=0mm, minimum size=.8mm, fill=red, shift={(2,0)}] (c5) at (g8) {};
					\node[circle, draw=red!80, inner sep=0mm, minimum size=.8mm, fill=red, shift={(2,0)}] (c9) at (g12) {};
					
					\node[circle, draw=red!80, inner sep=0mm, minimum size=.4mm, fill=red, shift={(0,-.5)}] (c6) at (c5) {};
					\foreach \i in {7,8}
					{
						\node[circle, draw=red!80, inner sep=0mm, minimum size=.4mm, fill=red, shift={(0,.25*(6-\i))}] (c\i) at (c6) {};
					}
					
					\node[circle, draw=red!80, inner sep=0mm, minimum size=.4mm, fill=red, label={right: $P_g$}] (c22) at (c7) {};

					\node[circle, draw=red!80, inner sep=0mm, minimum size=.8mm, fill=red, shift={(2,.5)}, label={right: $v_h$}] (c4) at (g1) {};
					\node[circle, draw=red!80, inner sep=0mm, minimum size=.8mm, fill=red, shift={(2,-.5)}, label={right: $v'_h$}] (c10) at (g19) {};
					
					\node[circle, draw=black!80, inner sep=0mm, minimum size=.4mm, fill=black, shift={(0,.5)}] (c3) at (c4) {};
					\foreach \i in {1,2}
					{
						\node[circle, draw=black!80, inner sep=0mm, minimum size=.4mm, fill=black, shift={(0,.25*(3-\i))}] (c\i) at (c3) {};
					}
					
					\node[circle, draw=black!80, inner sep=0mm, minimum size=.4mm, fill=black, shift={(0,-.5)}] (c11) at (c10) {};
					\foreach \i in {12,13}
					{
						\node[circle, draw=black!80, inner sep=0mm, minimum size=.4mm, fill=black, shift={(0,.25*(11-\i))}] (c\i) at (c11) {};
					}
					
					\draw  [line width=0.3mm, red] (c4) -- (c5) (c9) -- (c10) (g1) -- (c4) (g19) -- (c10);
					
					
					\foreach \i in {1,2,3,11,12,13}
					{
						\node[circle, draw=black!80, inner sep=0mm, minimum size=.4mm, fill=black, shift={(-8,0)}] (b\i) at (c\i) {};
					}
					
					\foreach \i in {6,7,8}
					{
						\node[circle, draw=blue!80, inner sep=0mm, minimum size=.4mm, fill=blue, shift={(-8,0)}] (b\i) at (c\i) {};
					}
					
					\node[circle, draw=black!80, inner sep=0mm, minimum size=.8mm, fill=black, shift={(-8,0)},label={left: $z_h$}] (b4) at (c4) {};
					\node[circle, draw=black!80, inner sep=0mm, minimum size=.8mm, fill=black, shift={(-8,0)},label={left: $z'_h$}] (b10) at (c10) {};

					\foreach \i in {5,9}
					{
						\node[circle, draw=blue!80, inner sep=0mm, minimum size=.8mm, fill=blue, shift={(-8,0)}] (b\i) at (c\i) {};
					}
					
					\node[circle, draw=blue!80, inner sep=0mm, minimum size=.4mm, fill=blue, label={right: $P_f$}] (b22) at (b7) {};

					\node[circle, draw=blue!80, inner sep=0mm, minimum size=.8mm, fill=blue, shift={(0,-1.2)}, label={left: $x'$}] (b14) at (b4) {};
					\node[circle, draw=blue!80, inner sep=0mm, minimum size=.8mm, fill=blue, shift={(0,1.2)},label={left: $y'$}] (b15) at (b10) {};
					
					\draw  [line width=0.3mm, black] (b4) -- (b14) (b15) -- (b10);
					\draw  [line width=0.3mm, blue] (b5) -- (b14) (b15) -- (b9);
					
					\draw  [line width=0.3mm, green] (b4) -- (g1) (b4) -- (c4) (b10) -- (g19) (b10) -- (c10);
					
					\draw  [line width=0.3mm, blue] (b14)  .. controls (4,8) .. (g7);
					\draw  [line width=0.3mm, blue] (b15)  .. controls (4,-6) .. (g13);

					\draw (8.5, 10) node[label=\textcolor{red}{$C'$}]  (Pc){};
					\draw (-6.5, 7.5) node[label=\textcolor{blue}{$C''$}]  (Pc){};

				\end{tikzpicture}
			}
			\caption{In red is cycle  $C'$ of Claim \ref{claim:towardsdrh2} and, in blue, cycle $C''$ of Claim~\ref{claim:pathsdh}.}
			\label{fig:F3}
		\end{figure}
		\begin{claim}\label{claim:towardsdrh2}
			Assume that we have $v\in R_{h,2}$. Let $u\in \dr(v)\cap P_d$ and assume that $Q_u:=u_{0}u_{1}u_{2}\dots u_{q}$ is a path witnesses that $u\in \dr(v)$ where $u_0=u$ and $u_q=v$. Then we have $u_2\in \{v_k,v'_k\}$.
		\end{claim}
		
		\proof
		By Claim~\ref{claim:towardsd} we have $u_1\in P_f$ and $u_2\in v_kP_hv_k'$. So it is enough to show that $u_2$ is not an internal vertex of $v_kP_hv_k'$ (and we can assume that this path has internal vertices). Clearly, any internal vertex of the path $v_kP_hv'_k$ are in the interior of the cycle $C'$ obtained from $C_{k,1}\cup C_{h,1}$ by deleting the internal vertices of $v_kP_hv'_k$ (see Figure~\ref{fig:F3}). On the other hand, since $f<h$ and $K_h$ is a component of $G_h=G-\cup_{j\le h-1}V(P_j)$ we have  $P_f\cap K_h=\varnothing$. As both the interior of $C_{k,1}$ and $C_{h,1}$ are in $K_h$, $P_f$ (and in particular $u_1$) is in the exterior of $C'$. As $u_2$ is the neighbour of $u_1\in P_f$, $u_2$ is not an internal vertex of $v_kP_hv'_k$.
		\hfill $\Diamond$

		\begin{claim}\label{claim:pathsdh}
			Assume that $v\in R_{h,2}$. If $|\dr(v)\cap P_d|\ge 6$, then $|\dr(v)\cap P_h|\le 7$.
		\end{claim}
		
		\proof Let $x,y\in \dr(v)\cap P_d$ be such that $d_{P_d}(x,y)$ is maximum, and $Q_x$ and $Q_y$ be paths witnessing that $x$ and $y$, respectively, belong to $\dr(v)$. Note that we have $3\le ||Q_x||, ||Q_y||\le 4$ and $Q_x, Q_y\subseteq K_d$. By Claim~\ref{claim:towardsdrh2}, if we take $Q_x=xx'x''\dots v$ and $Q_y=yy'y''\dots v$, then we have $x',y'\in P_f$ and $x'',y''\in \{v_k, v'_k\}$. First assume that $x''=y''$. As $P_d$ is isometric in $K_d$ and $xx'x''y'y$ is a walk of length 4 in $K_d$, we have that $||xP_dy||\le ||xx'x''y'y||\le 4$, which contradicts to our assumption that $|\dr(v)\cap P_d|\ge 6$. Thus $x''\neq y''$. We may assume that $x''=v_k$, $y''=v'_k$ and $v_k<_L v'_k$. We may also assume that $d_{P_h}(v_k,w_h)< d_{P_h}(v'_k,w_h)$. By the definition of $L$, we have $V(v_kP_hv'_k)<_L V(v'_kP_hw'_h-v'_k)$. Then it suffices to show that $|\dr(v)\cap w_hP_hv'_k|\le 7$ and $|\dr(v)\cap (v'_kP_hw'_h-v'_k)|=0$.
		
		Consider a vertex $z\in \dr(v)\cap w_hP_hv'_k$, and a path $Q_z$ witnessing $z\in \dr(v)$. We have $||Q_z||\le 4$ and $Q_z\subseteq K_h$. By Observation \ref{obs:startinki} we have $vQ_yv'_k-v'_k\subseteq K_k$, and by Observation~\ref{obs:subcomponent} we have $K_k\subset K_h$. Therefore $v'_k\in P_h\subseteq K_h$ and the walk $v'_kQ_yv$ lies in $K_h$. And this walk has length at most 2, so the walk $zQ_zvQ_yv'_k$ is of length at most 6 in $K_h$, and since $P_h$ is isometric in $K_h$ we obtain $||zP_hv'_k||\le ||zQ_zvQ_yv'_k|| \le 6$, therefore $|\dr(v)\cap w_hP_hv'_k|\le 7$. 
		
		Assume that $\dr(v)\cap (v'_kP_hw'_h-v'_k)\neq \varnothing$. Let $w\in \dr(v)\cap (v'_kP_hw'_h-v'_k)$ and $Q_w$ be a path witnessing $w\in \dr(v)$. By Observation~\ref{obs:intersectManagerForeman}, $Q_w$ intersects with either $v_kP_hv'_k$ or $z_kP_\ell z'_k$. As $V(v_kP_hv'_k)<_L V(v'_kP_hw'_h-v'_k)$, $Q_w\cap v_kP_hv'_k=\varnothing$, thus $Q_w\cap z_kP_\ell z'_k\neq \varnothing$. Let the cycle $C'':=v_kP_hv'_ky'P_fx'v_k$, then $w$ is in the exterior of $C''$. Next we will show that $z_kP_\ell z'_k$ is in the interior of $C''$. First, we know that $z_kP_\ell z'_k$ is adjacent to $P_k$ and by assumption $P_k$ is in $R_{h,2}$. Since $G$ is planar, $z_kP_\ell z'_k$ is also in $R_{h,2}$. The cycle $C_{h,2}$ together with two edges $v_kx'$ and $v'_ky'$ makes up three cycles, one of them is $C''$ which is the only one that contains both $v_k$ and $v'_k$. As $z_kP_\ell z'_k$ is adjacent to both $v_k$ and $v'_k$, $z_kP_\ell z'_k$ is in the interior of $C''$. Thus $Q_w$ intersects with $C''$. As $V(C'')\subseteq P_f\cup v_kP_hv'_k$ and $Q_w\cap v_kP_hv'_k=\varnothing$, we have that $Q_w\cap P_f\neq \varnothing$, which is impossible since $V(P_f)<_L V(P_h)$. Therefore, we have $\dr(v)\cap (v'_kP_hw'_h-v'_k)=\varnothing$.
		\hfill $\Diamond$

		If  we have $v\in R_{h,1}$, then by Lemma~\ref{lem:v in R_i1} we have $|\dr(v)\cap P_h|\le 8$, and by Claim~\ref{claim:pathdrh1} $|\dr(v)\cap P_d|\le 7$, obtaining \eqref{eq:dh}. Otherwise we use Claim~\ref{claim:pathd} and Claim~\ref{claim:pathsdh} to obtain \eqref{eq:dh}.
		
		The result now follows from \eqref{eq:Pk}, \eqref{eq:dh} and Claims~\ref{claim:ell}, \ref{claim:f}, \ref{claim:pathm}, and \ref{claim:pathn}. 
		\qed

We mention that the way to obtain the bound $|\dr(v)|\le 76$ is to prove that $|\dr(v)\cap(P_b\cup P_c\cup P_g)|\le 26$. The interested reader is referred to the first version of this paper~\cite{us} for details.

  \section{Concluding remarks: graphs of large treewidth}\label{sec:remarks}

Theorem \ref{thm:treewidth2bound7} gives an upper bound for $\chi(\stexact{G}{-2})$ when $G$ has treewidth at most 2. For graphs of larger treewidth, note that if $G$ has treewidth at most $t$, by using \eqref{eq:coltw} and Corollary \ref{coro:AM}, we obtain that $\chi(\stexact{G}{-2})\le   \col_2(G)\cdot2^{\col_2(G)-1} \le (t+1)\cdot 2^{t}$.  The fact that we need the term $2^t$ is attested by the following example. Take a clique on~$t$ vertices and an independent set $I$ on $2^t$ vertices. Make each vertex of $I$ adjacent to the whole clique, and for every pair of vertices in $I$ let them have a different signature on the edges that join them to the clique (and give any signature to the edges of the clique). Let $\hat{G}$ be the resulting signed graph. Then for each pair $u,v\in I$ there is a vertex $w$ in the clique such that the path $u,w,v$ is negative, and thus we have $\chi(\stexact{G}{-2})\ge 2^t$.

While this gives the right asymptotics, it would be nice to obtain tight bounds for $\chi(\stexact{G}{-2})$ when $G$ has bounded treewidth. As motivation we present the following short and direct proof of a weakening of Corollary \ref{coro:AM}.

 \begin{proposition}\label{thm:strongandcol}
		For every signed graph $\hat{G}$ we have
		$\chi(\stexact{G}{-2})\le \col_{2}(G)^2\cdot2^{\col_2(G)}$.
	\end{proposition}

 \begin{proof}
	Let $\sigma$ be the signature of $\hat{G}$, and $L$ be an ordering of $V(G)$ which witnesses that $\col_2(G,L)=\col_2(G)$. Moving along $L$, we give to each vertex $y$ a colour $a(y)\in [\col_2(G)]$, that is different from $a(x)$ for all $x\in \Reach_2 [G,L,y]\setminus\{y\}$, and a colour $b(y)\in [\col_2(G)^2]$, that is different from $b(z)$ for every $z$ in $$\mathrm{R}(y):=\{u \mid u\in \Reach_2[G,L,w] \mathrm{\,\, for\,\, some\,\, }  w\in \Reach_2[G,L,y] \}\setminus \{y\}. $$ To each vertex $y$ assign a vector $\alpha(y)$ with $\col_2(G)$ entries, such that the $i$-th entry is $+$ if $y$ has a neighbour $x<_Ly$ with $a(x)=i$ and $\sigma (xy)=+$, and is $-$ if $y$ has a neighbour $x<_Ly$ with $a(x)=i$ and $\sigma (xy)=-$; entries not yet defined are set to $+$. Note that $\alpha(y)$ is well defined, because if $x,z<_Ly$ are neighbours of $y$, then we either have $x\in \Reach_2[G,L,z]$ or $z\in \Reach_2[G,L,x]$, and thus $a(x)\ne a(z)$. Finally we assign to each $y\in V(G)$ a colour $c(y):=(\alpha(y),b(y))$, and claim that~$c$ is a proper colouring of $\stexact{G}{-2}$.
	
	Let $u,w\in V(G)$ be such that $d_G(u,w)=2$ and there is a negative $uw$-path of length 2 in $\hat{G}$. To prove our claim and finish the proof we now show that $c(u)\ne c(w)$. Let $u,v,w$ be a negative path in $\hat{G}$. If we have $v<_Lu$ and $v<_Lw$, then the $a(v)$-th entry of $\alpha(u)$ differs from the $a(v)$-th entry of $\alpha(w)$,   implying $c(u)\ne c(w)$. In every other case we either have $u\in\mathrm{R}(w)$ or $w\in\mathrm{R}(u)$, and thus $b(u)\ne b(w)$, implying $c(u)\ne c(w)$, as desired.  
\end{proof}

		\medskip
		{\bf Acknowledgment.} The project has received funding from the following grants. FONDECYT/ANID Iniciaci\'on en Investigaci\'on Grant 11201251, Programa Regional MATH-AMSUD \linebreak MATH230035 and MATH210008, program "Investissement d'Avenir" launched by the French Government and implemented by ANR, with the reference «ANR‐18‐IdEx‐0001» as part of its program «Emergence», Czech Science Foundation grant 25-16627S, NSFC (No. 12401471) and ZJNSFC (No. LQN25A010012).

	\end{document}